\documentclass[leqno]{amsart}
\usepackage{amssymb, amsmath, amsthm}
\usepackage[all]{xy}
\usepackage{color}

\newtheorem{thm}{Theorem}[section]

\newtheorem{cor}[thm]{Corollary}
\newtheorem{prop}[thm]{Proposition}
\newtheorem{lem}[thm]{Lemma}
\newtheorem{lemma}[thm]{Lemma}

\theoremstyle{definition} 
\newtheorem{defn}[thm]{Definition}
\newtheorem{eg}[thm]{Example}
\newtheorem{egs}[thm]{Examples}

\newtheorem{rmk}[thm]{Remark}

\newcommand{\zed}{\ensuremath{\mathbb Z}}

\newcommand{\Hom}{\mbox{\it Hom}}

\newcommand{\und}{\underline}

\newcommand{\into}{\hookrightarrow}
\newcommand{\iso}{\simeq}

\newcommand{\U}{\tilde{U}}
\newcommand{\V}{\tilde{V}}

\newcommand{\orbit}{\ensuremath{\mathcal{O}_G}}
\newcommand{\orbith}{\ensuremath{\mathcal{O}_H}}

\newcommand{\orbitx}{\ensuremath{\mathcal{O}_{G, X}}}

\begin{document}
\date{\today}

\title{Translation Groupoids and Orbifold Cohomology}
\author[D. Pronk and  L. Scull]{Dorette Pronk and Laura Scull}
\address{Department of Mathematics and Statistics,
Chase Building,
Dalhousie University, 
Halifax, NS, B3H 3J5, Canada}
\address{Department of Mathematics, Fort Lewis College,  1000 Rim Drive, 
Durango, Colorado 81301-3999    
USA
}
\email{pronk@mathstat.dal.ca, scull\_l@fortlewis.edu}

\thanks{Both authors are supported by NSERC discovery grants. Both authors thank 
the Fields Institute for its support and hospitality during
the Thematic Program on Geometric Applications of Homotopy Theory.
The first author also thanks Calvin College and Utrecht University  for their hospitality 
and the University of Chicago for its hospitality and financial support.}


\begin{abstract}
We show that the bicategory of (representable) orbifolds and good 
maps is equivalent to the bicategory of 
orbifold translation groupoids and generalized equivariant maps, giving a mechanism for transferring results from equivariant homotopy theory to the orbifold category.   As an application,  
we use this result to define orbifold versions of a couple of equivariant cohomology theories:   $K$-theory  and Bredon 
cohomology for certain coefficient diagrams.
\end{abstract}

\maketitle

\section{Introduction} 

Spaces with symmetries arise naturally in  many contexts, 
and have been studied from various points of view.    
Equivariant homotopy theory uses the tools of algebraic topology to study  
the category $G$-{\sf spaces}, consisting  of  
spaces with an action of the group $G$  and equivariant maps between them.    
Much of ordinary homotopy theory can be adapted and extended to this setting, 
although there are some important differences; see  \cite{alaska} 
for an overview of this theory.  From another point of view, there has been 
much recent interest in the study of orbifolds, which are something like manifolds 
but whose local structure is a quotient of an open subset of a Euclidean space  by a finite group 
action (\cite{sa56}; also  \cite{ALR, AR, LU}).  
Although many of the basic geometric structures are the same in both cases, the 
techniques of these two approaches have been rather different.  

The goal of this paper is to provide a way of moving between these points of view.  
One way to obtain an orbifold is to look at the action of a compact Lie group acting 
on a manifold with finite stabilizers.  In fact, a large class (perhaps all) of orbifolds 
can be described in this way \cite{HM}, although this description is not unique 
for a given orbifold.  Orbifolds that can be described this way are called 
representable.  We can try to  import equivariant invariants for these representable 
orbifolds.  In order to make this work, however, there are a couple of issues that 
need to be overcome.  The first is the fact that  the representation is not unique, 
and so in order to get invariants of the orbifold structure and not the particular 
representation, it needs to be checked that we get the same result for every 
representation.  The second, related, issue is that equivariant invariants are 
not  defined for non-equivariant descriptions of an  orbifold; and some orbifold maps 
may only be defined by  using  an alternate (potentially non-equivariant)  
description of the orbifold.  Thus we are faced with the possibility that a 
map between representable orbifolds may need to factor through an orbifold 
which does not come from a global group action, making it impossible to turn an 
equivariant invariant into a functor for the
orbifold category. 
 
In this paper, we prove that it is possible to represent every map
between representable orbifolds as an equivariant map, allowing us to define  equivariant invariants 
which are functorial for orbifold maps.  We also develop an explicit description 
of the non-uniqueness in the representation, making it practical to check which 
equivariant invariants will give orbifold invariants.
 
This non-uniqueness can be expressed in terms
of  Morita equivalences, generated by essential equivalences. 
These equivariant Morita equivalences are all compositions of 
certain specific forms of maps, and they  
satisfy the properties to allow us to form a bicategory of fractions
$$
\mbox{\sf Orbifolds}_{\mbox{\scriptsize\sf eqvar}}(W^{-1})
$$
where the Morita equivalences have become honest (internal) equivalences.
This same type of non-uniqueness is also
present in the description  of an orbifold in terms of an atlas of orbifold charts,
and  the category of orbifolds and good maps (or generalized maps)
is the bicategory of fractions 
$$
\mbox{\sf Orbifolds}_{\mbox{\scriptsize\sf atlas}}(W^{-1})
$$
of the category of orbifolds and atlas maps where the elements of the
class $W$ of essential equivalences have been `inverted'
to become equivalences.
We show that there is an equivalence
of bicategories for representable orbifolds,
$$
\mbox{\sf RepOrbifolds}_{\mbox{\scriptsize\sf eqvar}}(W^{-1})
\simeq\mbox{\sf RepOrbifolds}_{\mbox{\scriptsize\sf atlas}}(W^{-1}).
$$
Thus we have a more precise understanding of the relationship between the  
equivariant theory of the categories of $G$-spaces for various groups $G$,  
and the category of orbifolds;  this makes it possible to translate results between these settings,
and develop equivariant homotopy theory for orbifolds.

To demonstrate how this point of view can be applied,  we  show how the relationship between representable orbifolds 
and translation groupoids can be used to import $G$-equivariant cohomology theories to orbifolds; we discuss two such theories.  The first is topological  $K$-theory, defined using $G$-equivariant vector bundles; we show that this is in fact a well-defined cohomology theory on orbifolds.   
This  has been looked at in various other ways.  
Moerdijk  \cite{M02} has  shown that over the ring $\mathbb C$ of complex numbers,
this  can also be obtained as the equivariant sheaf cohomology of the inertia
groupoid $\Lambda(G)$ with values in the constant sheaf ${\mathbb{C}}$, and hence to prove that over $\mathbb{C}$, we get an orbifold cohomology theory.  This  approach could potentially be extended 
to  other coefficient systems by choosing the appropriate $\Lambda(G)$ sheaves, since there is a Leray spectral sequence relating the $K$-theory to the sheaf cohomology (via Bredon cohomology for certain coefficients).    In \cite{AR},  Adem and Ruan take an alternate approach and use   K-theory techniques to  get an orbifold invariant over the rationals $\mathbb Q$.   Our approach provides a more direct proof than either of these.

The second type of equivariant cohomology theories we consider are those defined by 
Bredon \cite{Br} with constant coefficients (coefficients which  do not depend on the space, 
only on the group $G$ and its orbit category).  
These Bredon cohomology theories are defined for coefficients given by diagrams of Abelian groups.    
We use our results to identify which of these coefficient diagrams actually give orbifold 
invariants, rather than depending on the equivariant representation used.  For these diagrams, 
we show that it is possible 
to define a notion of Bredon cohomology for representable orbifolds, 
depending only on the orbifold and not its equivariant presentation.  
Specifically, we describe a relation on these orbifold coefficient systems such that
if a $G$-space $X$ and an $H$-space $Y$ describe the same orbifold and 
$\underline{A}$ is a coefficient system on 
the orbit category of $G$, then there is a corresponding coefficient 
system on the orbit category of $H$ which gives the same cohomology groups. 

A related result was presented in Honkasalo's
paper \cite{Ho}. For a $G$-space $X$ with a coefficient
system $\und{A}$, Honkasalo constructs a sheaf $S(\und{A})$
on the orbit space $X/G$ such that the $G$-equivariant cohomology
on $X$ with coefficients in $\und{A}$ is isomorphic to the
sheaf cohomology of the orbit space $X/G$ with coefficients in $S(\und{A})$.
When applied to a representable orbifold, considered as a $G$-space,
it gives a relationship between the equivariant Bredon cohomology and the sheaf cohomology
of the underlying quotient space. This provides a nice alternative definition
of these cohomology groups. However, we should be careful not to read too much into
this description. For example, it does not imply that the Bredon cohomology is 
an invariant of the quotient space. The same topological space could be obtained 
as a different quotient $Y/H$ and there would not necessarily be an H-coefficient
system that would give rise to the same sheaf. For similar reasons, Honkasalo's
construction does not automatically give us an orbifold invariant (\cite{Ho}  does 
not  consider this question).  A sheaf which corresponds to a coefficient system 
for one representation does not need to correspond to a coefficient system in another representation,
as shown in Example \ref{counterex}.

Our approach gives a clearer idea of the relationship between 
the equivariant and orbifold phenomena, and is a blueprint for future 
applications of creating orbifold invariants out of equivariant ones.  
In a forthcoming paper we will construct an orbifold version of the
equivariant fundamental groupoid;  this is a category which 
has proved very useful in a variety of places in equivariant homotopy theory, 
including defining Bredon cohomology for twisted coefficients, 
obstruction theory  and studying equivariant orientations.  
We believe that this can be used to get analogous results for orbifolds, 
and perhaps lead to a characterization of the homotopy of the orbifold category. 

The organization of the paper is as follows.  Section \ref{S:back} gives an overview 
of the theory of orbifolds and how they are represented by groupoids.   
Section \ref{S:results}  gives the statements of our comparison results.   Section \ref{S:ktheory} gives the results on orbifold $K$-theory, and 
Section \ref{S:bredon}  contains the definitions of the Bredon cohomology for orbifolds.  
Sections \ref{S:main} and \ref{S:morepfs}  contain the deferred proofs 
of some of  the earlier results;   Section \ref{S:HS} contains supporting material for 
the proof of the main comparison theorem in Section \ref{S:main}.  

The authors thank Johann Leida for his stimulating conversations. Some of the questions that 
lead to this paper were inspired by his work on orbifold homotopy theory.
We also thank Ieke Moerdijk for his encouragement and for making them aware of 
some of the earlier literature related to this work.   Lastly, we thank Dev Sinha and the topologists at the University of Oregon for some helpful suggestions regarding equivariant $K$-theory.

\section{Background:  Orbifolds and Lie Groupoids} \label{S:back}

The classical definition of orbifolds (or V-manifolds) as first given by Satake,  
and developed  by Thurston and others,  is a  generalization of the definition of 
manifolds based on charts and atlases.  The difference is that the local neighbourhoods 
are homeomorphic to $U = \U/G$ where $G$ is a finite group acting on an open set 
$\U \subseteq {\mathbb R}^n$.
An orbifold can then be defined as a paracompact Hausdorff space $M$ together 
with an orbifold atlas, which is a locally compatible family of charts $(\U, G)$ such that the sets 
$\U/G$ give a cover of $M$.    The usual notion of equivalence of atlases through 
common refinement is used; details can be found in \cite{sa56,sa57}.
Note that the original definition required that all group actions be effective, 
but it  has been shown in recent papers (see for example, \cite{CR} or \cite{LU}) 
that it is often useful to drop this
requirement; we will not  require that $G$ acts effectively on $\U$. 

Working with orbifold atlases is cumbersome, particularly when dealing with maps 
between orbifolds.  Therefore an alternate way of representing orbifolds using groupoids 
has been developed.  It was shown in \cite{MP1} that every smooth orbifold can be represented 
by a Lie groupoid, which is determined up to essential equivalence.   
This way of representing orbifolds gives rise to a notion of orbifold map
which works well for homotopy theory \cite{MP1}. These maps have also 
been called `good' maps \cite{CR} or generalized maps.  This is the way we will 
approach the study of the orbifold category;  below, we review  some of the basic definitions.


\subsection{Lie Groupoids}\label{smgrpd}
A  groupoid is a (small) category in which all arrows are invertible.  
We think of the objects of the category as representing points in 
a geometric object, and the arrows as representing identifications.  
In order to reflect this information, we need to have a geometric 
structure present on our category.    
Therefore we work with Lie (or smooth) groupoids.

\begin{defn}
A (Hausdorff) {\em Lie groupoid} or {\em  smooth groupoid} $\mathcal{G}$ consists of smooth manifolds 
$G_0$ (the objects) and $G_1$ (the arrows)  together with the usual structure maps:  source and target  
$s, t\colon  G_1 \to G_0$, identity arrows determined by  $u\colon  G_0 \to G_1$, and composition 
$m\colon  G_1 \times_{s, G_0, t} G_1 \to G_1$, all given by smooth maps, such that $s$ (and therefore $t$)  
is a surjective submersion, and the usual diagrams commute (see, for example, Definition 4.1 in \cite{LU}).
\end{defn}

The following are examples of Lie groupoids:

\begin{egs}\label{gpdex}
\begin{enumerate}
\item  Any manifold can be viewed as a Lie groupoid by taking $G_1 = G_0 = M$, with only identity maps.   
\item  Any Lie group $G$ is a Lie groupoid with a single point $G_0 = *$, where composition of arrows is given by group multiplication.
\item
Let $G$ be a Lie group with a smooth  left action on a manifold $M$.
Then the translation groupoid $G\ltimes M$ is defined as follows.   
The objects are given by the manifold $M$ itself, and the arrows are defined by   $G\times M$.    
The source of an arrow $(g, x)$ is defined by $s(g, x) = x$, and the target by using the action of $G$ 
on $M$, $t(g, x) = gx$.    So $(g, x) $ is an arrow $x \to gx$.  
The other structure maps are defined by the unit  $u(x)=(e,x)$, where $e$ is the identity element in $G$,
and $(g',gx)\circ(g,x)=(g'g,x)$.
\end{enumerate}
\end{egs}

Now we define a category of Lie groupoids.   We use topologized versions of the usual category theory notions of functor and natural transformation; note that all maps are assumed to be smooth.

\begin{defn}
A {\em homomorphism} $\varphi\colon{\mathcal G}\rightarrow {\mathcal H}$ between Lie groupoids consists of a pair 
of maps $\varphi_0\colon G_0\rightarrow H_0$ and $\varphi_1\colon G_1\rightarrow H_1$, which commute with all the structure maps.

A {\em  natural transformation} or 2-cell between homomorphisms of Lie groupoids 
$\alpha\colon \varphi \Rightarrow \psi\colon {\mathcal G}\rightrightarrows{\mathcal H}$ 
consists of a map $\alpha\colon G_0\rightarrow H_1$ such that 
$s\circ\alpha=\varphi_0$, $t\circ\alpha=\psi_0$,
and $\alpha$ is natural in the sense that the following diagram commutes:
$$
\xymatrix{
G_1\ar[r]^-{(\psi_1,\alpha\circ s)} \ar[d]_{(\alpha\circ t,\varphi_1)} & H_1\times_{s,H_0,t}H_1\ar[d]^m
\\
H_1\times_{s,H_0,t}H_1\ar[r]_-m & H_1.
}
$$
\end{defn}

The category {\sf LieGpd} of Lie groupoids, homomorphisms, and natural transformations forms a 2-category.

\subsection{Essential Equivalences}\label{esseq}
We are thinking of  a groupoid as a representation of its underlying quotient space, 
encoding this space and its  singularity types.  However, this representation is not unique;  
the same quotient structure can be represented by different groupoids.  Therefore we need 
to introduce a notion of equivalence on the category of groupoids.

\begin{defn}
A homomorphism $\varphi\colon {\mathcal G}\rightarrow {\mathcal H}$ between Lie groupoids
is an {\em essential equivalence} when it satisfies the following two conditions.

1. It is {\em essentially surjective}, {\it i.e.}, the map 
\[ t\circ \pi_2\colon G_0\times_{H_0} H_1\rightarrow H_0\]   
from the manifold $G_0\times_{H_0} H_1 = \{ (x, h) \, | \, \phi_0 (x) = t(h) \}$ is a surjective submersion.

2. It is {\em fully faithful}, {\it i.e.}, the diagram 
$$\xymatrix@R=2em@C=2.5em{
G_1\ar[r]^-{\varphi_1}\ar[d]_{(s,t)} & H_1\ar[d]^{(s,t)}
\\
G_0\times G_0 \ar[r]_-{\varphi_0\times\varphi_0} & H_0\times H_0}
$$
is a pullback of manifolds.
\end{defn}

Thus an essential equivalence is a smooth equivalence of categories.  
Two groupoids ${\mathcal G}$ and ${\mathcal H}$ are {\em essentially equivalent} 
when there is a span of essential equivalences
$$
{\mathcal G}\leftarrow{\mathcal K}\rightarrow{\mathcal H}
$$
between them.  
In order to show that this is in fact an  equivalence relation,  
we use the notion of the (weak) fibre product of Lie groupoids. 
 
\begin{defn}     If $\phi\colon  \mathcal{H} \to \mathcal{G}$ and 
$\psi\colon  \mathcal{K} \to \mathcal{G}$ are homomorphisms of Lie groupoids, 
the {\em (weak) fibre product} $\mathcal{H} \times_{\mathcal{G}} \mathcal{K}$ (if it exists) is 
the following Lie groupoid.  The space of objects is the fibered product of manifolds
$H_0\times_{G_0}G_1\times_{G_0}K_0$. So an object is a  triple $(y, g, z)$ 
where $y \in {H}_0$, $z \in {K}_0$ and $g\colon  \phi(y) \to \psi(z)$ in $\mathcal{G}$.  
An arrow $(y, g, z) \to (y', g', z')$ consists of a pair $(h, k)$ of arrows
$h\colon  y \to y'$ in $\mathcal{H}$ and $k\colon z \to z'$ in $\mathcal{K}$ such that $g' \phi(h) = \psi(h) g$.
\end{defn} 

The fibre product introduced here has a 'weak' universal property of pullbacks for commuting diagrams 
of Lie groupoids and homomorphisms:  the square is  
only required to commute up to an (invertible) 2-cell.

Note that although source and target maps $s,t\colon G_1\rightrightarrows G_0$ 
are surjective submersions, this does not imply that $H_0\times_{G_0}G_1\times_{G_0}K_0$
is a manifold in general. The space $H_0\times_{G_0}G_1$ is a manifold,
but the map from this space into $G_0$ does not need to be transversal to the map
from $K_0$ into $G_0$. However, if at least one of the groupoid maps is an essential equivalence,
essential surjectivity gives that one of the maps involved in the last fibre product is again a submersion,
so we obtain another manifold.

It can also be shown that the fibre product  of an essential equivalence 
along any homomorphism is again an essential equivalence \cite{M88}; 
thus any zig-zag of essential equivalences may be shortened by taking a 
fibre product, and so by repeated shortening, replaced by a single span as above.

The class $W$ of essential equivalences between Lie groupoids also  satisfies the axioms needed  to form 
a bicategory in which  the essential equivalences have been inverted \cite{Pronkcomp}.  
In fact, the argument given for \'etale groupoids in \cite{Pronkcomp} works for Lie  groupoids as well.
So we can form the  bicategory of fractions $\mbox{\sf LieGpd}\,(W^{-1})$ as follows.  
The objects are the  Lie  groupoids as usual,
but a morphism ${\mathcal G}\rightarrow{\mathcal H}$ is a span of homomorphisms 
$$\xymatrix{{\mathcal G}&{\mathcal K}\ar[l]_\omega \ar[r]^\varphi & {\mathcal H}},$$
where $\omega$ is an essential equivalence. Such morphisms are also called {\em generalized maps}.  
Thus we are allowed to replace the source groupoid $\mathcal{G}$ with an essentially 
equivalent groupoid  $\mathcal{K}$ in defining our maps.  

 We define the composition of spans  using the fibre product construction.  
 In  showing that this fibre product gives a span of the right form, and so another generalized map,    
 the key point is again that the pullback of an essential equivalence along any homomorphism 
 is an essential equivalence.

A 2-cell $(\varphi,\omega)\Rightarrow(\varphi',\omega')$ in this bicategory is an equivalence class 
of diagrams of the form
$$\label{2-cells}
 \xymatrix@R=2em@C=3em{
	 &{\mathcal K}\ar[dl]_{\omega}\ar[dr]^\varphi
 \\
 {\mathcal G}\ar@{}[r]|{\alpha_1\Downarrow} 
	& {\mathcal L}\ar[u]_{\nu}\ar[d]^{\nu'}\ar@{}[r]|{\alpha_2\Downarrow} & {\mathcal H}
 \\
	 &{\mathcal K}'\ar[ul]^{\omega'}\ar[ur]_{\varphi'} &
 }
$$
 where $\omega\circ\nu$ and $\omega'\circ\nu'$ are essential equivalences.
 Note that since the essential equivalences satisfy a 2-for-3 property 
 (see Section \ref{S:morepfs},  Lemma \ref{2outof3} for a proof),  this is equivalent to 
 requiring that $\nu$ and $\nu'$ be essential equivalences.

Given an orbifold ${\mathcal M}$ with an orbifold atlas ${\mathcal U}$, we can define its groupoid representation ${\mathcal G}({\mathcal U})$ as follows. 
The space of objects is the disjoint union of the charts, 
$$G_0=\amalg_{\mathcal U}\U.$$
The space of arrows is a quotient of the space 
$$\coprod_{\begin{subarray}{c}
	\lambda_1\colon\U\into\V_1\\
	\lambda_2\colon\U\into\V_2\end{subarray}} \U,$$ 
where the disjoint union is
over pairs of atlas embeddings of $U$ into any charts.
The equivalence relation on the space of arrows is generated by a 
notion of local equivalence of pairs of embeddings;
the source and target maps on each copy of the charts are defined by the first and the second embedding respectively.
For further details, including  the definition of composition, 
the reader is referred to \cite{thesis}.

The local structure on these charts equips the resulting Lie groupoid with some special properties. 
In particular, a groupoid coming from an orbifold atlas  will satisfy the following conditions:

\begin{defn}  A Lie groupoid is:
\begin{enumerate}
\item
{\em\'etale} if its source map is a  local diffeomorphism:
\item
 {\em proper} if the map $(s,t)\colon G_1\rightarrow G_0\times G_0$ is a 
proper map ({\it i.e.}, it is closed with compact fibers).
\end{enumerate}
\end{defn}
  Note that if the source map is a local  diffeomorphism, this implies that the target map is also.

The notion of properness is preserved under the essential equivalence relation, but the notion of being \'etale
is not. 
This leads us to the following definition.

\begin{defn}
An {\em orbifold groupoid} is a groupoid which is essentially equivalent to a  proper \'etale Lie groupoid. 
\end{defn}

The construction outlined above shows that any orbifold can be represented by an orbifold groupoid.  
Conversely, given an orbifold groupoid $\mathcal{G}$, its orbit space can be given the structure of an orbifold \cite{MP1}.  
Therefore we consider such groupoids to be the orbifolds they represent.  
So {\sf Orbifolds} is the bicategory of orbifold groupoids with generalized maps as morphisms,  and 
equivalence classes of diagrams such as described above as 2-cells.

\section{Statement of Results:  \\ Representing Orbifolds by Translation Groupoids} \label{S:results}

In order to make a bridge between orbifolds and equivariant homotopy theory, 
we are interested in representing orbifolds by a particular type of  Lie groupoid:  
the translation groupoids $G \ltimes M$ coming from the action of a Lie group $G$ 
on a manifold $M$, described in Example \ref{gpdex}, part (3).  
It turns out that many, possibly  all, orbifolds can be represented this way.
Satake showed  that every effective orbifold can be obtained as a quotient of a manifold
by the action of a (not necessarily finite) compact Lie group acting with finite isotropy \cite{sa57}.
Unfortunately, his proof does not go through for non-effective orbifolds. 
However, a partial result was obtained by  Henriques and Metzler \cite{HM}; 
their Corollary 5.6 shows that  all 
orbifolds for which all the ineffective isotropy groups have trivial centers are representable.   
It is conjectured that all orbifolds are representable, but this has not been proven.  

For the remainder of this paper, we restrict our attention to those orbifolds 
that are representable, so that we can work with their translation groupoids.    
In this section, we give the statements of our results showing that we can form 
a bicategory of representable orbifolds using only translation groupoids and equivariant maps.    
The proofs of these statements are generally constructive and sometimes rather long, 
and so we will defer many of them until  Sections \ref{S:main} and  \ref{S:morepfs}. 

\begin{defn}
An {\em equivariant map} $G\ltimes X\rightarrow H\ltimes Y$ between translation groupoids consists of a pair $(\varphi,f)$, 
where $\varphi\colon G\rightarrow H$
is a group homomorphism and $f\colon X\rightarrow Y$ is a $\varphi$-equivariant smooth map, 
{\it i.e.}, $f(gx)=\varphi(g) f(x)$ 
for $g\in G$ and $x\in X$.
\end{defn}  

We will denote the 2-category of smooth translation groupoids and equivariant maps
by {\sf EqTrGpd}.  

In order to represent orbifolds by objects in this category, we  
need to identify essentially equivalent groupoids as before.  
So we want to invert the essential equivalences in  {\sf EqTrGpd}, 
and show that we can form a  bicategory by defining maps using spans 
as in the previous section.  Again, the key to making this process work
is the fact that we can form the fibre product of translation groupoids 
and get another translation groupoid, in such a way  that the pullback of 
an essential equivalence along an equivariant map is another essential equivalence.   

\begin{lem} \label{L:pb}
In a fibre  product of Lie groupoids
$$
\xymatrix{
{\mathcal P}\ar[r]^\zeta\ar[d]_\xi \ar@{}[dr]|{\cong} & G\ltimes X\ar[d]^{\psi}
\\
H\ltimes Y\ar[r]_{\varphi} &{\mathcal K},
}
$$
the groupoid  ${\mathcal P}$ is again a translation groupoid.  
Moreover, its structure group is $G\times H$, and $\zeta$ and $\xi$ 
are equivariant maps, where the group homomorphisms are the appropriate projections.
\end{lem}

The proof examines the explicit construction of $\mathcal{P}$ to verify the claims about it, and is given in Section \ref{S:morepfs}.

\begin{cor}  For every pair of equivariant maps 
$$\xymatrix{H\ltimes Y\ar[r]^{(\varphi,f)} & K\ltimes Z & \ar[l]_{(\psi,w)}G\ltimes X}$$
where $(\psi,w)$ is an essential equivalence, there is a commutative square
$$
\xymatrix{
L\ltimes P\ar[r]^{(\zeta,z)} \ar[d]_{(\xi,v)} & G\ltimes X\ar[d]^{(\psi,w)}
\\
H\ltimes Y\ar[r]_{(\varphi,f)} & K\ltimes Z 
}
$$
where  $(\xi, v)$ is an essential equivalence.  
\end{cor}

\begin{proof}  
We combine the result from  \cite{M88} that the fibre product of an essential
equivalence is again an essential equivalence with Lemma \ref{L:pb}, to show that 
the usual  pullback essential equivalence  is again in the 
2-category {\sf EqTrGpd}.
\end{proof}

Thus, we can again replace any zig-zag of essential equivalences with a single span; 
and also define the composition of spans, which is unitary and associative up 
to coherent isomorphisms. Finally, it is straightforward to adjust 
the proof that the class $W$ of general essential equivalences
satisfies the conditions to admit a bicategory of fractions $\mbox{\sf LieGpd}(W^{-1})$,
to show that the class of equivariant essential equivalences, which we will  again call $W$,  gives rise 
to a well-defined bicategory $\mbox{\sf EqTrGpd}(W^{-1})$.

Now we want to show that for representable orbifolds, restricting to the equivariant maps 
of the category  $\mbox{\sf EqTrGpd}\,(W^{-1})$ does not lose important information; 
that is, $\mbox{\sf EqTrGpd}\,(W^{-1})$ is equivalent to the bicategory 
 $\mbox{\sf LieGpd}_{\mbox{\scriptsize\sf Tr}}(W^{-1})$,
the full sub-bicategory of $\mbox{\sf LieGpd}\,(W^{-1})$ on translation groupoids.

\begin{thm} \label{T:main}
The inclusion functor $\mbox{\sf EqTrGpd}\hookrightarrow \mbox{\sf LieGpd}_{\mbox{\scriptsize\sf Tr}}$ 
induces an equivalence of bicategories
$$
\mbox{\sf EqTrGpd}(W^{-1})\simeq\mbox{\sf LieGpd}_{\mbox{\scriptsize\sf Tr}}(W^{-1}),
$$
when restricted to orbifold groupoids.
\end{thm}

The proof involves replacing  generalized maps and 2-cells by  equivariant ones between translation groupoids, 
in such a way that this induces the desired equivalence of bicategories.    
For instance, for  a generalized map $G \ltimes X \leftarrow \mathcal{K} \to H \ltimes Y$ 
between orbifold translation groupoids, we construct an isomorphic span of equivariant maps  
$$
\xymatrix{
G\ltimes X &\ar[l]_{(\psi,w)} K\ltimes Z\ar[r]^{(\varphi,f)} & H\ltimes Y.
}
$$  
To  construct $K \ltimes Z$ (and also the corresponding replacements for the 2-cells), 
we  make use of an alternate way of describing maps between orbifold Lie groupoids, 
based on groupoid bundles,  developed by Hilsum and Skandalis \cite{HS}.  
Some background on Hilsum-Skandalis maps is given in Section \ref{S:HS}, 
and Section \ref{S:main} gives the proof of Theorem \ref{T:main}.

Thus we can work with just the equivariant maps between translation groupoids, 
with the equivariant essential equivalences inverted.    
There are a couple of obvious forms of equivariant maps which are essential equivalences:  
if we have a $G$-space $X$ such that a normal subgroup $K$ of $G$ acts freely on $X$, 
then it is easy to see that the quotient map 
\begin{equation}\label{quotientform}
G\ltimes X \rightarrow G/K\ltimes X/K,
\end{equation}
is an essential equivalence.    Similarly, for any (not necessarily normal) subgroup $K$ of a group 
$H$ and $K$-space $Z$, we can induce up to get an $H$-space 
$H \times _K Z = G \times Z / \sim$, where $[hk, z] \sim [h, kz]$ for any $k \in K$.  
Then the inclusion $Z \to H \times_K Z$ defined by $z \to [e, z]$ gives an essential equivalence 
\begin{equation}\label{inclform}
K\ltimes Z\rightarrow H\ltimes (H\times_K Z),
\end{equation}

It turns out that these are the only forms of equivariant weak equivalences we need to deal with,
since they generate all other equivariant essential equivalences.

\begin{prop}    \label{L:form}
Any equivariant essential equivalence is a composite of maps of the forms (\ref{quotientform}) 
and (\ref{inclform}) described above.
\end{prop}

We defer the proof until  Section \ref{S:morepfs}.
 
Thus we have an explicit description for the weak equivalences in $\mbox{\sf EqTrGpd}(W^{-1})$.

\section{Orbifold K-Theory}  \label{S:ktheory}

One example of an equivariant cohomology theory 
that has been extensively studied is equivariant K-theory;  
see \cite{se, alaska} for an introduction to this theory.   
This is defined geometrically using $G$-equivariant vector bundles 
for compact Lie groups $G$, and has many applications.    
Elementary properties of these equivariant vector bundles can 
be combined with our results on representation of orbifolds to give an 
easy proof that $K$-theory is actually an orbifold invariant. 
This result has been  proved by Adem and Ruan in \cite{AR} over the rationals and Moerdijk \cite{M02} 
over the complex numbers by various other approaches, as discussed in the introduction.

  We briefly recall the definition of equivariant $K$-theory for a compact Lie group $G$.    
A $G$-vector bundle over a $G$-space $X$ is a vector bundle $\xi:  E \to X$ such that the 
total space $E$ also has  a $G$-action making the projection map an equivariant map, 
and such that $G$ acts linearly on fibres $\xi^{-1}(x) \to \xi^{-1}(gx)$.    
For a compact space $X$, the equivariant $K$-theory $K_G(X)$ is defined as the 
Grothendieck group of finite dimensional $G$-vector bundles over $X$.
Tensor product makes this into a ring.
  
We can extend this  to a cohomology theory on $G$-spaces as follows.    
We can  define a reduced version of the $K$-theory group for spaces 
with a $G$-fixed basepoint by $\tilde{K}_G(X) = \ker [(K_G(X)  \to K_G(*)]$  
(for unbased spaces, we adjoin a disjoint fixed basepoint);   then equivariant 
Bott periodicity holds for $\tilde{K}_G$:  with complex coefficients, 
$\tilde{K}_G(X_{+}) \iso \tilde{K}_G(S^V \wedge X_+)$ for a complex representation $V$;  
similarly for an $8$-dimensional real representation.   
Thus we  can  define a  cohomology theory by  $\tilde{K}^n_G(X) = \tilde{K}_G(\Sigma^n X)$ \cite{se, alaska}.

\begin{prop}  \label{P:ktheory} 
Suppose $\mathcal{X}$ is an orbifold represented by a translation groupoid $G \ltimes X$.  
Then the  equivariant $K$-theory  group  $K_{G}(X)$  
is independent of the  representation.
\end{prop}

\begin{proof}
By Theorem \ref{T:main} and Proposition \ref{L:form},  
it is sufficient to check that the definition is invariant under the two forms of  
change-of-group essential equivalences  (\ref{quotientform}) and (\ref{inclform}).  

The geometric definition of the group $K_G(X) $ makes it easy to see 
that if $X$ is a free $G$-space, then the $G$-vector bundles correspond 
to ordinary vector bundles over the quotient space $X/G$.  More generally,   
if $H $ acts freely on $X$, then $K_G(X) \iso K_{G/H}(X/H)$.  
Therefore this is invariant under quotient maps $X \to X/H$ for free $H$-actions (\ref{quotientform}).  
Similarly, $K_G(G \times_HY) \iso K_H(Y)$, since any $G$-bundle over $G\times_HY$ is determined 
by its underlying $H$-bundle over $Y$.   
Therefore this definition is also independent of the inclusion change-of-groups  (\ref{inclform}).

Thus the group  $K^0_G$ is an orbifold invariant.   
For the general cohomology theory, we need only observe  that $S^1 \wedge (G \times_H X)_+ $ 
is canonically isomorphic to \\  $G\times_H (S^1 \wedge X_+)$, 
and hence the higher $K$-groups are also invariant under this equivariant change-of-groups.
\end{proof}

\section{Orbifold Bredon Cohomology}  \label{S:bredon}

In this section, we use the results on representing orbifolds via equivariant spaces  
to develop a definition of Bredon cohomology for orbifolds.  
Throughout, we will again assume that all groups are  compact Lie groups, 
and that all subgroups are closed.
 
Bredon cohomology takes its inspiration from the idea that we should view a $G$-space 
as being described by the diagram of its fixed points 
$\{ X^H \, | \, hx = x, \, \forall \, h \in H \}$ for the various subgroups $H$ of $G$.    
The natural inclusions and $G$-action give morphisms between these sets.  
These can be organized by the  {\em orbit category } $\orbit$, 
which has the canonical $G$-orbit types $G/H$  as its objects, 
with  all equivariant maps between them.   
These equivariant maps can be described concretely as composites 
of maps of the form $G/H \to G/\alpha H \alpha^{-1}$ defined 
by $gH \to g\alpha H$, and projection maps $G/H \to G/H'$ for $H \subseteq H'$.  
Since the fixed set $X^H$ can also be described as the equivariant mapping space 
$\mbox{\em Hom}_G(G/H, X)$ from the canonical orbit $G/H$, 
we immediately see that the fixed sets form a (contravariant) 
functor to the category of spaces   
$\Phi X \colon  \orbit^{\mbox{\scriptsize\it op}} \to {\sf Spaces}$ defined  by   $\Phi X(G/H) = X^H$.   

Many of the usual algebraic invariants of spaces can then be composed with the functor $\Phi X$ to give diagrams 
of Abelian groups indexed by $\orbit$.  
Moreover, any homotopy invariants will result in diagrams indexed by the 
homotopy category $h \orbit$, which has homotopy classes of 
equivariant maps as its morphisms.    
Thus the home for many equivariant invariants is the category 
${\sf Ab}^{h\orbit^{op}}$ of functors $h\orbit^{\mbox{\scriptsize\it op}} \to {\sf Ab}$, which we call  
{\em coefficient systems}.

The category of coefficient systems can be used to define a cohomology theory as follows.  
We have a chain complex of coefficient systems  
$\und{C}_*(X)$ defined by  \[ \und{C}_n(X) (G/H)  = C_n(X^H / WH_0)\]  
where $WH_0$ is the identity component of the Weyl group $NH/H$.    
Then for any coefficient system $\und{A}$, the maps between diagrams $\und{C}_n(X)$ 
and $\und{A}$ are defined by the natural transformations, and these form an Abelian group  
\[ {C}^n_{h\orbit} (X; \und{A}) = \Hom_{h\orbit}(\und{C}_n(X), \und{A}).\]    
Together these give a graded Abelian group $C^*_{h\orbit}(X; \und{A})$.  
The boundaries on the chains in the fixed point spaces $\und{C}_*(X)$ induce a differential on this,  
and so we obtain a cochain complex $C^*_{h\orbit}(X; \und{A}) $.  
The Bredon cohomology of $X$ is then the cohomology of this  complex:   
\[ H^*_{G}(X; \und{A}) = H^*_{h\orbit}(X; \und{A})= H^*(C^*_{h\orbit}(X; \und{A}) )\] 
and this construction satisfies the axioms for an equivariant cohomology theory on $G$-{\sf Spaces} \cite{Br}.  

We observe that for a given $G$-space $X$, although the Bredon cohomology is 
defined on diagrams indexed by all closed subgroups of $G$, it really  
only depends on  isotropy subgroups of $X$.

 \begin{prop}  \label{P:res} 
Let $h\orbitx$ be the full subcategory of $h\orbit$ on objects $G/H$ such 
that $H$ is an isotropy group of $X$.   
Then $H^*_{h\orbit} (X, \und{A}) = H^*_{h\orbitx}(X, r_X \und{A})$ 
where $r_X\und{A}$ is the restriction of the diagram $\und{A}$  to  $h\orbitx$.
 \end{prop}

\begin{proof} 
The chain complex  $\und{C}_*(X)$ is generated by cells of an equivariant 
$G$-CW decomposition of cells of $X$; such cells are of the form $G/H \times D^n$ 
for some isotropy subgroup $H$, and contribute a summand $\und{G/H} \times \zed $ 
to $\und{C}_n(X)$, where the diagram $\und{G/H}$ is defined by 
$\und{G/H} (G/L) = \pi_0(G/H)^L$.  But $\pi_0(G/H)^L = \Hom_{h\orbit}(G/L, G/H)$ 
and so a Yoneda argument shows that 
$\Hom_{h \orbit}(\und{G/H} \times \zed, \und{A}) = \Hom_{Ab}(\zed, \und{A}(G/H))$.    So  
$\mbox{\it Hom}_{h\orbit} (\und{C}_*(X), \und{A}) \iso \coprod_{H} \mbox{\it Hom}\,(\zed, \und{A}(G/H))$ 
where $H$ runs through the isotropy groups of $G$.   
Because only isotropy groups are involved, this is exactly the same as 
$\mbox{\it Hom}_{h\orbitx} (\und{C}_*(X), r_X \und{A})$.   
(See \cite{W}).
\end{proof}

\begin{cor}\label{rX}
If $r_X\und{A} = r_X\und{B}$ then $H^*_{h\orbit}(X, \und{A}) = H^*_{h\orbit}(X, \und{B})$.
\end{cor}

Alternatively, this also follows from a result by Honkasalo \cite{Ho}, which 
describes the Bredon cohomology of a $G$-space $X$ with coefficient system
$\und{A}$ as the sheaf cohomology of the orbit space $X/G$ with coefficients
in the sheaf $S(\und{A})$ with stalks $S(\und{A})_{\overline{x}}\cong \und{A}(G/G_x)$. 
Honkasalo's result suggests that for representable orbifolds, it may be possible to define 
Bredon cohomology as an orbifold invariant.    
However, it does not completely identify which coefficient systems will give an orbifold invariant
as we will see below.

Since representable orbifolds can be represented as translation groupoids, 
we can apply  the above definitions to a particular translation 
groupoid representation  of an orbifold, and obtain cohomology groups.  
In order to be a true orbifold invariant, however, we need  to ensure 
that these groups do not depend on the representation, {\it i.e.}, 
that the definition of orbifold Bredon cohomology sees 
only structure associated to the orbifold, and not the particular translation groupoid representing it.  
The results of the previous section tell us exactly what is required 
to be an invariant of the orbifold:  
we need a definition that is invariant under the change-of-group 
essential equivalences described in Proposition \ref{L:form}.  
In particular, if $G\ltimes X$ and $H\ltimes Y$ are two essentially equivalent orbifold translation 
groupoids, we need a correspondence between coefficient systems on $h{\mathcal O}_G$ and $h{\mathcal O}_H$,
which will give the same cohomology groups for $X$ and $Y$ respectively.

We will see that this is not possible for all coefficient systems, because 
some coefficient systems give invariants that are not orbifold invariants.      
One way that the equivariant theory may distinguish finer information than 
that carried by the orbifold structure is to differentiate between 
(disjoint) fixed point sets which have isomorphic isotropy and in fact are part of the
fixed point set of the same subgroup in some representations, as in the following example.

\begin{eg}\label{counterex}
Let $Q$ be the orbifold represented as the quotient of the circle $S^1$ by the action of
$D_2 = {\mathbb Z}/2\oplus{\mathbb Z}/2=\langle\sigma_1\rangle\oplus\langle\sigma_2\rangle$, 
where $\sigma_1$ and $\sigma_2$ act by horizontal and vertical reflection respectively.  
The points with non-trivial isotropy groups are the north and south poles and 
the east and west poles, each with isotropy group ${\mathbb Z}/2$. Then the 
subgroup generated by $\sigma_1\sigma_2$ acts freely, so we can take the quotient
to obtain a new representation for $Q$ as $S^1/\langle\sigma_1\sigma_2\rangle\cong S^1$ with an
action of $D_2/\langle\sigma_1\sigma_2\rangle\cong {\mathbb Z}/2$.
In this case ${\mathbb Z}/2$ acts by reflection with two fixed points. In this second presentation,
the subgroup fixing these points is the same, where there pre-images in the first presentation
had distinct isotropy subgroups (which were isomorphic, but not the same, or even conjugate,
as subgroups of $D_2$).  So an orbifold coefficient system cannot attach distinct 
Abelian groups to these subgroups.
\end{eg}

Given an orbifold
$\mathcal X$ represented by a $G$-space $X$ 
and a coefficient system
$$\underline{A}\colon h{\mathcal O}_G\rightarrow {\sf Ab},$$ we want to construct coefficient systems
$\Psi_*\underline{A}\colon h{\mathcal O}_H\rightarrow {\sf Ab}$
and $\Phi^*\underline{A}\colon h{\mathcal O}_K\rightarrow {\sf Ab}$
for all essential equivalences
$\Psi\colon G\ltimes X\rightarrow H\ltimes Y$
and $\Phi\colon K\ltimes Z\rightarrow G\ltimes X$, in such a way that
the essential equivalences induce
isomorphisms between the cohomology groups with coefficients in 
the corresponding coefficient systems.

The example above shows that this is not always possible.
However, we will show that the issue of taking a quotient
by a freely acting subgroup, which was the cause of the problems
in this example, is the only
one we need to address. Moreover, it is always possible to determine
from the given representation whether this issue will arise.
So we can give a characterization of orbifold coefficient systems
which only depends on the given representation.

In general, the previous section shows that if $K$ is a normal 
subgroup of $G$ acting freely on $X$, 
then the $G$-space $X$ is orbifold equivalent to the $(G/K)$-space $X/K$.
Therefore we have to be careful when  $X$ has fixed sets $X^L$ and $X^{L'}$  
associated to subgroups $L$ and $L'$ whose  projections in $G/K$ are  the same.    
In this case, in the quotient space $X/K$ these become part of the same fixed set
$(X/K)^{LK/K}$, and 
so an orbifold cohomology theory must treat these the same.  Looked at another way, 
we must be able to deduce all the information contained in 
the orbifold Bredon cohomology groups with coefficient in a system $\und{A}$ 
on $h \orbit$ from the Bredon cohomology of the quotient $X/K$  
defined with diagrams on  $h \mathcal{O}_{G/K}$.   

Some of this happens automatically, as we observe from the following.

\begin{lemma}\label{unique}
Let $X$ be a $G$-space, and $K$  a normal subgroup of $G$ which acts freely on $X$.
If some point $x \in X$ is fixed by two subgroups $L$ and $L' $ of $G$  with the same projection in $G/K$,  
then   $L = L'$. 
 \end{lemma}
 
\begin{proof}  
Since $L$ and $L'$ have the same projection in $G/K$, then for any $\ell' \in L'$ we must have $\ell' = \ell k$ for some $k \in K$.   
Suppose that  $\ell x = x$ and $\ell' x = x$; so $\ell k x = x$.  Then  $x = \ell^{-1}x$ and thus  $kx = \ell^{-1}x = x$; 
so $k$ must also fix $x$.  Since $K$ acts freely, we conclude that $k = e$ and so $L = L'$.
\end{proof} 

Thus, if $x \in X$ is a lift of $\overline{x} \in X/K$ which is fixed by $\overline{L} \subseteq G/K$, 
there is a unique subgroup $L \subseteq G_x$ lifting $\overline{L}$;
the equivariant Bredon cohomology does not have a chance to distinguish between 
different lifts of $\overline{L}$ at the point $x$, since there is a unique lift $L$
such that $x\in X^L$.

In order to behave as an orbifold invariant,  we also need the Bredon cohomology  
to treat all lifts $x$ of $\overline{x}$ equally; again this follows from elementary group theory.
Of course, if $L$ is an isotropy subgroup of a point $x$, then the conjugates $gLg^{-1}$ 
are isotropy subgroups of the points $gx$ in the orbit; since conjugation is an isomorphism in $h \orbit$, 
the values of any coefficient system $\und{A}$ are isomorphic at all conjugates:  
$\und{A}(G/L) \cong \und{A}(G/gLg^{-1})$.  Moreover, we have the following  result.

\begin{lemma}  \label{L:triv} Suppose that  $K, H$ are subgroups of $G$ such that $K$ is normal and $K \cap H = \{ e \}$.  
If the conjugation action by $K$ fixes  $H$, then in fact  $K$  acts trivially on $H$.
\end{lemma} 

\begin{proof}
Suppose $k \in K \cap NH$; so $khk^{-1} = \hat{h}$.  Then $kh = \hat{h}k$ so $kh\hat{h}^{-1} = \hat{h}k\hat{h}^{-1}$.  
But $K$ is normal so $\hat{h}k\hat{h}^{-1} = \hat{k}$.  Then $kh\hat{h}^{-1} = \hat{k}$ and so $h\hat{h}^{-1} = k^{-1}\hat{k}$ 
is in $K \cap H = \{e\}$.  So $h = \hat{h}$.
\end{proof}
 
Thus, if $K $ acts freely on a $G$-space $X$, and $\overline{x}\in (X/K)^{\overline L}$, then for any lifts 
$x$ and $x'$ of $\overline{x}$ with (uniquely specified) subgroups $L\le G_x$ and $L'\le G_{x'}$ respectively,
lifting $\overline{L}$,  the points $x$ and $x'$ will differ by some $k \in K$, with $x' = kx$; so $L' = kLk^{-1}$.  
This conjugation by $k$ induces an isomorphism between $\und{A}(G/L)$ and 
$\und{A}(G/kLk^{-1})=\und{A}(G/L')$ which does not depend on the choice of $k$ by Lemma \ref{L:triv}.
So we can consider this to be a canonical identification.  
Thus we have  $\und{A}(G/L) = \und{A}(G/L')$ in this case, and any Bredon cohomology will treat these the same.

We conclude that many of the necessary identifications for an orbifold invariant are already present in any coefficient system.
It is possible, however, to have two {\em non-conjugate} isotropy subgroups 
$L$ and $L'$ which project to the same 
subgroup of $G/K$; thus we do need to place a restriction on our diagrams.

We will say that a subgroup $L\le G$ is in the {\em isotropy lineage} of $X$ if it is a subgroup of an 
isotropy group $G_x$ for some point $x\in X$, or equivalently, if $X^L\neq\emptyset$.

\begin{defn} 
We say that a diagram $\und{A}$ is an {\em orbifold coefficient system } 
if it satisfies the following: If $K$ is a normal subgroup of $G$ acting freely on 
$X$, then for any subgroups $L$ and $L'$ in the isotropy lineage of $X$ with $LK /K = L'K/K$ in 
$G/K$, we  have $\und{A}(G/L) = \und{A}(G/L')$;  similarly, any two 
structure maps of $h\orbit$ between isotropy lineage subgroups which project to the same 
structure map in $h \mathcal{O}_{G/K}$ must be identical in the
coefficient system.
\end{defn}

This identifies which diagrams have a chance to define an orbifold Bredon cohomology.   
The condition only becomes a real restriction when there are non-conjugate subgroups 
in the isotropy lineage with the same projection.

We will now show that our definition of orbifold Bredon cohomology is independent of the translation groupoid 
representation used.  Thus, for any two translation groupoids with a 
change-of-groups essential equivalence
between them,  we will identify which coefficient system on the 
one orbit category corresponds to a given orbifold coefficient system on the other.

In fact, Moerdijk and Svensson \cite{MS} have considered the issue of change-of-groups 
maps for Bredon cohomology as a special case of a more general study of changing 
the underlying diagram shape in diagram cohomology.    
If $\phi\colon G \to K$ is any group homomorphism, 
then $\phi$ induces a map $h\orbit \to h\mathcal{O}_K$  
defined on objects by $G/H \mapsto K/\phi(H)$.  This in turn induces a map of coefficient systems 
$\phi^* \colon {\sf Ab} ^{h\mathcal{O}_K^{op}} \to {\sf Ab} ^{h\orbit^{op}}$ 
defined by pre-composition with $\phi$.      (The authors of \cite{MS} state their 
results for discrete groups and use the orbit categories $\mathcal{O}_G$, but their arguments
extend to compact Lie groups when one uses the homotopy orbit categories $h\orbit$.) We will use  the following:  

\begin{prop}[\cite{MS}, Proposition 1.2] \label{P:MS} 
If $\phi\colon  G \to K$ is any group homomorphism and $X$ is a $G$-CW complex, 
then \[ H^*_K(K \times_{\phi, G} X, \und{A}) \iso H^*_G(X, \phi^*\und{A}) \]  
where $K \times_{\phi,G} X = K \times G / (k, gx) \sim (k \phi(g), x)$.
\end{prop}

The two particular group homomorphisms we are interested in are those coming from the 
essential equivalence change-of-group maps of Proposition \ref{L:form}.  The first form is a  projection 
$\pi:  G \to G/K$ for a normal subgroup $K \subseteq G$ which acts freely on the space $X$.  
In this case, $G/K \times_{\phi, G} X \cong X/K$ and 
$\pi^*\und{A}(G/L) = \und{A} ((G/K)/(L/K \cap L)) = \und{A}( (G/K) / (LK/K)) $.  

The second form is the inclusion  $i:  H \into G$  
for any  subgroup $H \subseteq G$, where $G \times_{i,H} X = G \times_H X$ 
is the usual space induced by the extension of groups.  In this case,   
$i^*\und{A}(H/L) = \und{A}(G/L)$; 
thus, $i^*$ just restricts the diagram to the subgroups contained in $H$.

Motivated by these observations and Corollary \ref{rX}, we make the following definition. 
 
\begin{defn}  \label{D:orbicoeff}
Suppose we have an orbifold ${\mathcal X}$ represented by a translation groupoid ${\mathcal G}=G\ltimes X$.
We define an equivalence relation on orbifold coefficient systems, denoted {\em orbifold equivalence},  
generated by the following.
\begin{itemize}  
\item  
If $K$ is a normal subgroup of $G$ which acts freely on $X$,  an orbifold system    
$\und{A}$ on $h \mathcal{O}_{G/K}$ is equivalent to $\pi^*\und{A}$  on $h\orbit$ 
where $$\pi^*\und{A}(G/H) = \und{A}((G/K)/(HK/K)).$$
\item   
If $H \subseteq G$ is any subgroup, then an orbifold system
 $\und{A}$ on $h\orbit$ is equivalent to $i^* \und{A}$ on $h\orbith$ where 
$$i^*\und{A}(H/L) = \und{A}(G/L).$$
\item
Two orbifold coefficient systems $\und{A}$ and $\und{B}$ on $h{\mathcal O}_G$ are equivalent if
$$r_X\und{A} = r_X\und{B}.$$  
\end{itemize} 
 If $[\und{A}]$ is an equivalence class of orbifold coefficient systems 
represented by $\und{A}$ on $h\orbit$, then the {\em Bredon cohomology} $H^*_{Br}(\mathcal{X}, [\und{A}])$ 
is defined by $H^*_{\orbit}(X, \und{A})$.  
\end{defn}

Note that for any orbifold coefficient system $\und{A}$, the induced systems $\pi^*\und{A}$ and $i^*\und{A}$ 
are again orbifold coefficient systems.

It is clear from the definition of the equivalence relation that 
for any essential equivalence $G\ltimes X\rightarrow L\ltimes Y$, and any orbifold  coefficient system 
on $h{\mathcal O}_L$, there is an equivalent system on $h{\mathcal O}_G$.
The following lemmas give us the other direction, namely that for any 
orbifold coefficient system on $h{\mathcal O}_G$ there is an equivalent system on $h{\mathcal O}_L$.

\begin{lemma}\label{subgroup}
Suppose we have an orbifold $\mathcal{X}$ represented by a translation groupoid 
$\mathcal{H} = H \ltimes X$, and let $\und{A}$ be an orbifold coefficient system on $h\orbith$.
For any larger group $G$ containing $H$ as subgroup, there is an orbifold coefficient system $\und{B}$ 
for the $G$-space $G\times_HX$
on $h\orbit$ such that $i^*\und{B}$ is equivalent to $\und{A}$.
\end{lemma}

\begin{proof}
Let $i\colon H\hookrightarrow G$ be the inclusion of groups.
Define the coefficient system $\und{B}:=i_*(\und{A})$ on $h\orbit$ 
in the following way. Let $L$ be a subgroup of $G$.
If $L$ is not in the isotropy lineage of $G$, we define $\und{B}(G/L)=0$.
If on the other hand, $(G \times_H X)^L\neq\emptyset$, let $[g,x]$ be a point in this fixed point set.
In that case $g^{-1}Lg$ is a subgroup of $H$, since it keeps the point $[e,x]$ fixed.
So we define $\und{B}(H/L):=\und{A}(G/(g^{-1}Lg))$.
It is not hard to see that $\und{B}$ defined this way is an orbifold
coefficient system when $\und{A}$ is and that $r_X i^*\und{B} = r_X \und{A}$.
\end{proof}

\begin{lemma}\label{quotientgroup}
Suppose we have an orbifold $\mathcal{X}$ represented by a translation groupoid 
$\mathcal{G} = G \ltimes X$, and that $K$ is a normal subgroup of $G$ which acts freely on $X$. 
For every orbifold coefficient system
$\und{A}$ on $h{\mathcal O}_{G}$ there is an orbifold  coefficient system  $\und{B}$ on $h{\mathcal O}_{G/K}$
such that $\pi^*\und{B}$ is equivalent to $\und{A}$.
\end{lemma}

\begin{proof}
Given the system $\und{A}$ on $h\orbit$, define the system $\und{B}:=\pi_*(\und{A})$ on $h{\mathcal O}_{G/K}$ 
as follows. Given a subgroup $\overline{L}$ of $G/K$, if $\overline{L}$ is not in the isotropy
lineage of $X/K$, then define $\und{B}((G/K)/\overline{L})=0$.  Otherwise, 
choose a  point $\overline{x}_{\overline{L}}\in (X/K)^{\overline{L}}$
and a point $x_{\overline{L}}\in X$ such that $\pi(x_{\overline L})=\overline{x}_{\overline{L}}$.
Let $L'\subseteq G$ be the unique subgroup of the isotropy group of $x_{\overline{L}}$
such that $L'K/K=\overline{L}$.
Define $\und{B}((G/K)/\overline{L})=\und{A}(G/L')$.    
Note that the choice of $L'$ may depend on 
the choice of $\overline{x}$, and up to conjugacy with an element in $K$ on the choice of $x_{\overline L}$, 
but the value of $\und{B}((G/K)/\overline{L})$ does not, 
because $\und{A}$ is an orbifold coefficient system.

Now we need to define structure maps for the coefficient system $\und{B}$ for the non-trivial portion of the diagram.
If we have two subgroups $\overline{L}_1 \subseteq \overline{L}_2$ in the isotropy lineage of $X/K$, 
we know that $L_2'$ has some subgroup $L_1''$ which projects to $\overline{L}_1$; and then 
$\und{B}((G/K)/\overline{L}_1) = \und{A}(G/L_1'')$, since this value does not depend on which lift is chosen.  
Then we can define the structure map associated to the projection map 
$(G/K)/(\overline{L}_1)\rightarrow(G/K)/(\overline{L}_2)$ to be the structure map 
$\und{A}(G/L_2) \to \und{A}(G/L_1'')$. 

For any left multiplication map $(G/K)/\overline{L}\rightarrow (G/K)/(\overline{g}(\overline{L})\overline{g}^{-1})$
in $h{\mathcal O}_{G/K}$, we need to define a morphism
$\und{B}((G/K)/(\overline{g}(\overline{L})\overline{g}^{-1}))\rightarrow\und{B}((G/K)/\overline{L})$.
Note that $\overline{L}$ is in the isotropy lineage of $X/K$ if and only if  ${\overline{g}(\overline{L})\overline{g}^{-1}}$ is.  
When both are in the isotropy lineage, let 
$L'$ be the chosen lift of $\overline{L}$, and $L''$  the chosen lift of $\overline{g}\overline{L}\overline{g}^{-1}$.    
We need a morphism
\begin{equation}\label{conj1}
\und{A}(G/L'')=\und{B}((G/K)/(\overline{g}(\overline{L})\overline{g}^{-1}))\rightarrow\und{B}((G/K)/(\overline{L}))=\und{A}(G/L').
\end{equation}
Now if we pick any pre-image $g$ of $\overline{g}$, then $gL'g^{-1}$ also projects to $\overline{g}\overline{L}\overline{g}^{-1}$ in $G/K$; so  $\und{A}(G/(gL'g^{-1})')=\und{A}(G/L'')$, 
since $\und{A}$ is an orbifold coefficient system. 
So the structure map in (\ref{conj1}) is defined as the structure map induced by left multiplication with $g$ in the orbit category,
\begin{equation}\label{conj2}
\und{A}(G/L'')=\und{A}(gL'g^{-1})\rightarrow\und{A}(G/L').
\end{equation}
Note that this map does not depend on the choice of $g$ such that $gK=\overline{g}$, since structure maps related to 
multiplication with elements of $k$ correspond to the canonical identifications of the groups in the diagram.

Thus, we have defined a coefficient system $\und{B}$ on $h\mathcal{O}_{G/K}$.  We complete this proof by showing that $r_X\pi^*\und{B}=r_X\und{A}$.
For an isotropy group $H\le G$, we have $\pi^*\und{B}(G/H)=
\und{B}((G/K)/(HK/K))=\und{A}(G/(HK)')$, where $(HK)'K/K=HK/K$ and 
$(HK)'$ is an isotropy group. Since $\und{A}$ is an orbifold coefficient system,
this implies that $\und{A}(G/(HK)')=\und{A}(G/H)$, so $\pi^*\und{B}(G/H)=\und{A}(G/H)$.
\end{proof}


Note that if $G \iso H \oplus K$, there are two canonical ways to obtain an equivalent system 
of coefficients on $h\orbit$ from one on $h\orbith$:  
using  $i_*$ for $i:  H \into G$ or using $\pi^*$ for $\pi:  G \to G/K \iso H$.  The result is the same either way,    
since  $LK / K \iso L / L\cap K$, and $L = (L\cap H)(L \cap K)$ so 
$L / (L \cap K) = (L \cap H) (L \cap K)/(L \cap K) \iso L \cap H / (L \cap K \cap H) = L \cap H$.
 
\begin{prop}  \label{P:br}
For any orbifold system of coefficients $\und{A}$,
$H^*_{Br}(\mathcal{X}, [\und{A}])$ is well-defined; that is, 
it does not depend on what translation groupoid is used to represent $\mathcal{X}$.  
 \end{prop}
 
\begin{proof}  By Theorem \ref{T:main} and Proposition \ref{L:form},  
it is sufficient to check that the definition is invariant under the two forms of  
change-of-group essential equivalences  (\ref{quotientform}) and (\ref{inclform}).  
Let $\varphi \colon G\ltimes X\rightarrow H\ltimes Y$ be such an essential equivalence.
By the definition of our equivalence relation and Lemmas \ref{subgroup} and
\ref{quotientgroup}, for any coefficient system $\und{A}$ on $h\orbit$ there is a coefficient system 
$\und{B}$ on $h\orbith$ which is equivalent to $\und{A}$, and conversely, for any coefficient system 
$\und{B}$ on $h\orbith$ there is an equivalent system $\und{A}$ on $h\orbit$.   So it is sufficient
to check that any of the pairs of systems that generate the equivalence relation give isomorphic
cohomology groups. 

Both of the first two cases of the equivalence relation follow directly from Proposition \ref{P:MS}.
The third case follows from Corollary \ref{rX}.
\end{proof}

Thus we have a way of defining orbifold Bredon cohomology under mild
restrictions on the coefficient systems. These restrictions are needed
because in some sense we have taken the limit of the Bredon cohomologies
for all the different equivarant representations of a given orbifold.
If one would like to take all coefficient systems for all representations
of the orbifold into account, one would need to consider a kind of colimit
construction. We plan to address these issues in more detail in a future
paper.

\begin{eg}
An example of a coefficient system which satisfies the conditions necessary to give rise to an orbifold cohomology theory is  $\und{R}_G$ defined by the representation rings:   
such a coefficient diagram  is defined by  $\und{R}_G(G/H)= R(H)$ on $h\orbit$.   
The structure maps of this diagram are induced by the conjugate $G$-action on itself:    
for a map $\alpha\colon  G/H \to G/(\alpha H \alpha^{-1})$, and a representation $V$ 
of $\alpha H \alpha^{-1}$, we simply pre-compose with the conjugation to get a representation of $H$.  
Similarly, for a projection $G/H \to G/H'$ for $H \subseteq H'$ and a representation $W$ of $H'$, 
we can restrict to an action of $H$ via the inclusion.   
$\und{R}_G$ is clearly an orbifold system of coefficients, since the value at $G/H$ only depends on $H$.
\end{eg}
 
Moreover, up to isomorphism of cohomology groups, 
it doesn't matter which translation groupoid we start with to represent 
our orbifold, as shown by the following proposition.

 \begin{prop}  \label{P:rep} 
Suppose $\mathcal{X}$ is an orbifold.
Then if $\mathcal{X}$ is  represented by two different translation groupoids 
$G \ltimes X$ and $H \ltimes Y$, then the orbifold coefficient systems 
$\und{R}_G$ and $\und{R}_H$ are orbifold equivalent and give rise to 
isomorphic Bredon cohomology
groups.
\end{prop}

\begin{proof}

  By Theorem \ref{T:main}, Proposition \ref{L:form} and Proposition \ref{P:br}, it is sufficient to check 
for any  change-of-group essential equivalences   of the forms (\ref{quotientform}) and (\ref{inclform}),
that the representation coefficient system on the domain is orbifold equivalent to the representation coefficient system on
the codomain.     

Let $i\colon  H \into G$ induce the essential equivalence $H \ltimes X \to G \ltimes (G \times_H X)$. 
Then $i^*\und{R}_{G} $ gives a diagram on $h\orbith$ by  
restricting to the subgroups  $K$ contained in $H$, and similarly restricting to those 
structure maps induced by the action of $H$; that is,  the restriction is exactly  $i^*\und{R}_G = \und{R}_H$.

Let $K \subseteq G$ act freely on $X$, inducing the essential equivalence $G \ltimes X \to G/K \ltimes X/K$.
In this case,   $\pi^*(\und{R}_{G/K})$ is {\em not} isomorphic to $ \und{R}_G$.  
However, we will now show that  $r_X \pi^*\und{R}_{G/K} = r_X \und{R}_G$.   

Because  $K$ acts freely on $X$, we know that if $H$ is an isotropy subgroup of $X$,  then $H \cap K = \{ e \}$.    
Therefore $\pi^*\und{R}_{G/K}(G/H) = \und{R}_{G/K}((G/K) / (HK/K)) = R(HK/K) = R(H/K\cap H) = R(H)$;
so the entries of $r_X \pi^*\und{R}_{G/H} $ and $ r_X \und{R}_G$ agree.
We also need to show that the structure 
maps of these two diagrams agree.

If $H \subseteq H'$ are two isotropy subgroups of $X$ in $G$, then the structure map of the 
projection $G/H \to G/H'$ is defined by restricting the $H'$ actions of the representations in $R(H')$ to $H$.
In $\pi^*\und{R}_{G/K}$, the structure map of the projection $G/H \to G/H'$ is induced by 
considering the projection $(G/K) / (HK / K) \to (G/K) / (H'K/K)$, and so comes from restricting the 
$H'K /K$ action to the subgroup $HK/K$.  But again, this is just isomorphic to the inclusion $H \subseteq H'$.  
So these structure maps are the same on the representation rings.

The diagram  $\pi^*\und{R}_{G/K}$ also has structure maps induced on the representations by 
the conjugation action of $G/K$ on its subgroups. In particular, any conjugation action of an 
element of $K$ is trivial in $\und{R}_{G/K}$ and therefore also in the induced diagram $\pi^*\und{R}_{G/K}$.
On the other hand, the diagram  $\und{R}_G$ has potentially more of these conjugation actions, 
coming from the action  of the larger group $G$ on its subgroups.   However, these extra morphisms 
are actually trivial:   any  $k \in K$ which sends a subgroup $H$ to itself, acts trivially on $H$ by 
Lemma \ref{L:triv}, and so any such structure map on $R(H)$ is already trivial.  Similarly, if $k$ takes $R(H)$ to   
the isomorphic ring $R(kHk^{-1})$, all such elements $k\in K$ must give the same isomorphism.
Thus all morphisms in the conjugation action of $G$ on the isotropy subgroups of $X$ factor through 
$G \to G/K$, and so the structure maps and thus the diagrams  $r_X \pi^*\und{R}_{G/K}$ 
and $r_X \und{R}_G$ are equal as desired.

Note that a similar argument can be applied to show that for any inclusion of groups $i\colon  H \into G$
and orbifold groupoid $H\ltimes Y$,
$r_{G\times_HY}i_*\underline{R}_H=r_{G\times_HY}\underline{R}_G$. And similarly, that for any
essential equivalence of the form $\pi\colon G \ltimes X \to G/K \ltimes X/K$,
$r_{X/K}\pi_*\underline{R}_G=r_{X/K}\underline{R}_{G/K}$.

So if both $G\ltimes X$ and $H\ltimes Y$ represent the same orbifold, let $\tilde{\underline{R}}_H$ be the 
coefficient system on $h\orbit$ obtained by moving $\underline{R}_H$ along some zig-zag of essential equivalences connecting 
$G\ltimes X$ and $H\ltimes Y$. From the argument above we derive that $r_X \und{R}_G=r_X \und{\tilde{R}}_H$.
So $H_{h\orbit}^*(X,\und{R}_G)=H_{h\orbit}^*(X,\und{\tilde{R}}_H)=H_{h\orbith}^*(Y, \und{R}_H)$.
\end{proof}

This particular Bredon cohomology theory is of interest because these are the ground coefficients for equivariant K-theory. K-theory is not itself a Bredon cohomology theory,  since its value evaluated at a point is not concentrated in degree zero.   However, K-theory is related to this particular  Bredon cohomology via an equivariant Atiyah-Hirzebruch spectral sequence.

Note that the techniques of this paper do not necessarily guarantee that there is a canonical
isomorphism between the cohomology groups in Proposition \ref{P:rep}.
The issue is that it is not clear whether two parallel essential equivalences of orbifold groupoids give rise to the same isomorphism, even if they give rise to the same maps between
the corresponding quotient spaces. (If there is an invertible 2-cell between them, this is the case, but for noneffective orbifolds it is not clear whether such a 2-cell needs to exist.)

A possible approach to this question would involve Honkasalo's description of 
these cohomology groups in terms of the sheaf cohomology of the quotient space.
A complete proof would require a construction of Honkasalo's sheaf $S(\und{R})$ 
based on the isotropy groups alone, without any reference to a representation
$G\ltimes X$ of the orbifold, together with a canonical isomorphism 
$S(\und{R}_G)\cong S(\und{R})$ for any such representation. This would require some careful
arguments about chosen embeddings of atlas charts, and fall outside the scope of the present paper.

\section{Proofs I:  Background on Hilsum-Skandalis Maps} \label{S:HS} 

The remainder of this paper consists of the deferred proofs of the results already discussed. 
We begin with supporting material for the proof of Theorem \ref{T:main}.

In using Lie groupoids  to represent geometric objects like orbifolds, often one ignores the bicategory structure and instead considers
the category $[\mbox{\sf LieGpd}]$ of Lie  groupoids with isomorphism classes of homomorphisms, 
and its corresponding category of fractions $[\mbox{\sf LieGpd}][W^{-1}]$ with respect to isomorphism classes of essential equivalences.
The advantage of considering this category rather than its 2-categorical refinement is that there is a nice description of the
morphisms in terms of groupoid bundles.
The resulting morphisms are called Hilsum-Skandalis maps \cite{Hae,Pr}. 
In \cite{MM}, Moerdijk and Mr\v{c}un give a description
of the correspondence between isomorphism classes of generalized maps and Hilsum-Skandalis maps, 
which we will use in the proof of Theorem \ref{T:main}.  
In this section, we summarize these constructions in order to fix our notation.

\begin{defn}
A {\em left ${\mathcal G}$-bundle} over a manifold $M$ is a manifold $R$ with smooth maps
$$\xymatrix@C=1.5em@R=1.5em{ 
R \ar[r]^\rho \ar[d]_r &M\\G_0}$$ and a left ${\mathcal G}$-action $\mu$ on $R$, with anchor map 
$r\colon R\rightarrow G_0$, such that  $\rho(gx) = \rho(x)$ for any $x\in R$
and any $g\in G_1$ with $r(x) = s(g)$. 

Such a bundle $R$ is {\em principal} if
\begin{enumerate}
\item 
$\rho$ is a surjective submersion, and 
\item
the map $(\pi_1, \mu) \colon R \times_{G_0} G_1\rightarrow R \times_M R$, sending $(x, g)$ to $(x, gx)$, is a
diffeomorphism.
\end{enumerate}
\end{defn}

A {\em Hilsum-Skandalis map} ${\mathcal G}\rightarrow {\mathcal H}$ 
is represented by a principal right ${\mathcal H}$-bundle $R$ over $G_0$
$$
\xymatrix@R=2em@C=2em{
 R \ar[r]^\rho \ar[d]_r &G_0\\H_0
}
$$
which also has  a left ${\mathcal G}$-action (along $\rho$), which commutes with the ${\mathcal H}$-action.
So we have that $$r(gx)=r(x),\, \,  \rho(xh)=\rho(x),\mbox{ and }g(xh)=(gx)h,$$ for any $x\in R$, $g\in G_1$ and $h\in H_1$
with $s(g)=\rho(x)$ and $t(h)=r(x)$. 
Moreover, since the $\mathcal{H}$-bundle is  principal, 
$\rho$ is a surjective submersion, and the map $R \times_{H_0} H_1\rightarrow R \times_{G_0} R$
is a diffeomorphism. We denote this map by $(R,\rho,r)\colon{\mathcal G}\rightarrow{\mathcal H}$.
Two principal right ${\mathcal H}$-bundles with left ${\mathcal G}$-action 
represent the same Hilsum-Skandalis map if and only if
they are diffeomorphic as ${\mathcal H}$- and ${\mathcal G}$-bundles.

{\em Composition} of Hilsum-Skandalis maps is defined by a tensor product construction over the middle groupoid.
Let $(R,\rho,r)\colon {\mathcal G}\rightarrow{\mathcal K}$ and 
$(Q,\theta,q)\colon{\mathcal K}\rightarrow {\mathcal H}$ be two Hilsum-Skandalis maps.
Then the space $\theta\circ\pi_2=\rho\circ\pi_1\colon R\times_{K_0}Q\rightarrow K_0$ 
has a right ${\mathcal K}$-action, defined by
$(x,y)k=(xk,k^{-1}y)$, for $k\in K_1$, $x\in R$, and $y\in Q$, with $\theta(y)=t(k)=\rho(x)$.
Denote the orbit space of this action by $R\otimes_{\mathcal K}Q$.  
Then we define the composition  
$$
(Q,\theta,q)\circ(R,\rho,r)=(R\otimes_{\mathcal K} Q,\rho\circ\pi_1,q\circ\pi_2).
$$

\begin{eg} The  left ${\mathcal G}$-bundle 
$$
\xymatrix@R=1.5em@C=1.5em{
G_1\ar[r]^s\ar[d]_t & G_0
\\
G_0
}
$$
defined by composition of morphisms is principal, and has also a right ${\mathcal G}$-action with anchor map  $s$ (again, by composition).
We denote this bundle by $U({\mathcal G})$.
The bundles of the form $U({\mathcal G})$ represent  identity morphisms in the sense that if $(R,\rho,r)$ is a  Hilsum-Skandalis map $\mathcal{H} \to \mathcal{G}$, then 
 $U({\mathcal H})\circ(R,\rho,r)\cong(R,\rho,r)\cong(R,\rho,r)\circ U({\mathcal G}).$
\end{eg}

\begin{defn}
A Hilsum-Skandalis map  $(R,\rho,r)$ is a {\em Morita equivalence} 
when it is both a principal ${\mathcal G}$-bundle and a principal ${\mathcal H}$-bundle.
\end{defn}

We can translate between Hilsum-Skandalis maps and our homomorphisms of Lie groupoids as follows.  
Let $\varphi\colon{\mathcal G}\rightarrow{\mathcal H}$ be a homomorphism. 
Then let $R_\varphi=\varphi^*U({\mathcal H})=G_0\times_{\varphi,H_0,t}H_1$.
This space has the following smooth functions to $G_0$ and $H_0$: 
\begin{equation}\label{HSmap}
\xymatrix@C=3em{G_0 & \ar[l]_-{\pi_1}R_\varphi\ar[r]^-{s\circ\pi_2} & H_0,}
\end{equation}
where $\pi_1$ and $\pi_2$ are the projection maps.
Moreover, it is a principal right ${\mathcal H}$-, and left ${\mathcal G}$-bundle with the following actions: 
$$g(x,h)h'=(t(g), \varphi_1(g)hh')$$ for $x\in G_0$, $g\in G_1$ and $h',h\in H_1$, with $s(g)=x\mbox{ and }t(h')=s(h).$
So (\ref{HSmap}) denotes a Hilsum-Skandalis map 
$(R_\varphi,\pi_1,s\circ\pi_2)\colon{\mathcal G}\rightarrow{\mathcal H}$.

Conversely, a Hilsum-Skandalis map 
$(R,\rho,r)\colon {\mathcal G}\rightarrow{\mathcal H}$ gives rise to a generalized map:
$$
\xymatrix@C=2em{
{\mathcal G}& \ar[l]_-{\tilde\rho} ({\mathcal G}\times {\mathcal H})\ltimes{R}\ar[r]^-{\tilde r} & {\mathcal H}
}
$$
where 
$$
(({\mathcal G}\times{\mathcal H})\ltimes{R})_0=R, \mbox{ and }
(({\mathcal G}\times{\mathcal H})\ltimes{R})_1=G_1\times_{s,G_0,\rho}R\times_{r,H_0,s}H_1,
$$ 
with $s(g,x,h)=x$, $t(g,x,h)=gxh^{-1}$, and
$m((g',gxh^{-1},h'),(g,x,h))=(g'g,x,h'h)$.

The homomorphisms $\tilde\rho$ an $\tilde{r}$ are defined by 
$$\tilde{\rho}_0(x)=\rho(x),\quad \tilde{\rho}_1(g,x,h)=g$$ and $$\tilde{r}_0(x)=r(x), \quad\tilde{r}_1(g,x,h)=h.$$

These constructions satisfy the following properties.

\begin{thm} 
{\em \cite{Mr}} The homomorphism $\tilde{r}$ is an essential equivalence if and only if $(R,\rho,r)$ is a Morita equivalence.
\end{thm}

\begin{thm} {\em \cite{Mr}}
The category of Lie groupoids with Hilsum-Skandalis maps forms a category of fractions  
for the category of Lie groupoids
with equivalence classes of homomorphisms relative to the essential equivalences.
\end{thm}

\section{Proofs II:  Proof of Theorem \ref{T:main}} \label{S:main}

We want to show that the bicategory of orbifold translation groupoids and equivariant 
maps in $\mbox{\sf EqTrGpd}(W^{-1})$ is equivalent to the full sub-bicategory of ${\sf LieGpd}(W^{-1})$ 
on representable orbifold groupoids; so we need to show that we can restrict to equivariant maps.

Let 
\begin{equation}\label{span}
\xymatrix{G\ltimes X&\ar[l]_-\upsilon {\mathcal K} \ar[r]^-\varphi  & H\ltimes Y}
\end{equation}
be a generalized map between translation groupoids. 
The fact that $\upsilon$ is an essential equivalence does not imply that
$\mathcal K$ is a translation groupoid. However, we will  show that it is isomorphic in ${\sf LieGpd}(W^{-1})$
to a generalized map of the form
\begin{equation}\label{span'}
\xymatrix{G\ltimes X&\ar[l]_-{\omega} L\ltimes Z\ar[r]^-\psi  & H\ltimes Y},
\end{equation}
where $\omega$ is a smooth equivariant essential equivalence 
and $\psi$ is a smooth equivariant map.
We will use the Hilsum-Skandalis representation of generalized 
maps as described in the previous section to construct 
the generalized map in (\ref{span'}).

\begin{prop}\label{morphisms}
Let ${\mathcal G}=G\ltimes X$ and ${\mathcal H}=H\ltimes Y$ be orbifold translation groupoids.
Any generalized map 
$$\xymatrix@1{{\mathcal G}&\ar[l]_\upsilon {\mathcal K} \ar[r]^\varphi  & {\mathcal H}}$$ 
is isomorphic in the
bicategory $\mbox{\sf LieGpd}[W^{-1}]$ to a generalized map of the form 
$\xymatrix@1{{\mathcal G}&\ar[l]_\omega {\mathcal L} \ar[r]^\psi  & {\mathcal H}}$ 
where ${\mathcal L}$ is a translation groupoid
and both $\omega$ and $\psi$ are equivariant maps. 
Moreover, ${\mathcal L}$ may be chosen such that its structure group is $G\times H$ and the 
group homomorphisms of $\omega$ and $\psi$ are the appropriate  projections onto $G$ and $H$.
\end{prop}

\begin{proof}
Let $R_\upsilon$ and $R_\varphi$ be the principal bundles 
corresponding to the homomorphisms $\upsilon$ and $\varphi$ respectively, 
as in (\ref{HSmap}) in Section \ref{S:HS}.
So 
$$
R_\upsilon=\upsilon^*(U{\mathcal G})=K_0\times_{X}(G\times X)
$$ 
and its elements can be represented as triples
$(z,g,x)$ with $z\in K_0$, $g\in G$, and $x\in X$, such that $\upsilon_0(z)=gx$.
Note that given $z$ and $g$, we have that $x=g^{-1}\upsilon_0(z)$,
so $$R_\upsilon\cong K_0\times G.$$
The projection map $\pi_1\colon R_\upsilon\rightarrow K_0$ is a surjective submersion, since it is the pullback
of the target map  $t\colon G\times X\rightarrow X$, $t(g,x)=gx$, which is a surjective submersion.
The anchor maps for the bundle structures on $R_\upsilon\cong K_0\times G$ are now
$$\xymatrix{K_0&\ar[l]_-{\pi_1} R_\upsilon\ar[r]^-{r_\upsilon} & X,}$$
where $r_\upsilon(z,g)=g^{-1}\upsilon_0(z)$.
The right $\mathcal G$-action and left $\mathcal K$-action are defined by
$$k\cdot(z,g)\cdot (g',g'^{-1}g^{-1}\upsilon_0(z))=(t(k),\pi_1\upsilon_1(k)gg').$$

Since $\upsilon$ is an essential equivalence, $R_\upsilon$ is also a principal ${\mathcal G}$-bundle,
representing a Hilsum-Skandalis map ${\mathcal G}\rightarrow{\mathcal K}$. 
As such we will denote it by $R_\upsilon^{-1}$; 
the space is the same, but the actions are reversed. (Recall that a left (resp. right) action can be turned into
a right (resp. left) action by acting by the inverses of the elements.)

The principal ${\mathcal K}$-bundle $R_\varphi$ is defined analogously.  We consider the composition of the two Hilsum-Skandalis maps represented by $R_\upsilon^{-1}$
and $R_\varphi$.
The principal bundle for the composition is obtained as a quotient of the pullback 
$$Q=R_\upsilon^{-1}\times_{K_0}R_\varphi\cong G\times K_0\times H.$$
The right ${\mathcal K}$-action on the projection map $Q\rightarrow K_0$ is 
defined by $$(g,z,h)\cdot k=(\pi_1\upsilon_1(k^{-1})g,\, s(k),\, \pi_1\varphi_1(k^{-1})h),$$ for $k\in K_1$ with
$t(k)=z$.
The quotient of $Q$ by this action is $R_\upsilon^{-1}\otimes_{\mathcal K}R_\varphi$.
This space has the following bundle maps into $X$ and $Y$:
$$\xymatrix{X & \ar[l]_-{q_{\upsilon}}\ar[r]^-{q_\varphi}R_\upsilon\otimes_{\mathcal K}R_\varphi & Y,}$$
where 
$$
q_\upsilon(g,z,h)=g^{-1}\upsilon_0(z)\mbox{ and }  q_\varphi(g,z,h)=h^{-1}\varphi_0(z).
$$
These maps are well-defined on equivalence classes, since 
\begin{eqnarray*}q_\upsilon(\pi_1\upsilon_1(k^{-1})g,s(k),\pi_1\varphi_1^{-1}(k^{-1})h)&=&
[\pi_1\upsilon_1(k^{-1})g]^{-1}\upsilon_0(s(k))\\
&=&g^{-1}\pi_1\upsilon_1(k^{-1})^{-1}\upsilon_0(s(k))\\
&=&g^{-1}\pi_1(\upsilon_1(k))\upsilon_0(s(k))\\
&=&g^{-1}\upsilon_0(t(k))\\
&=& g^{-1}(\upsilon_0(z))\\
&=&q_\upsilon(g,z,h).
\end{eqnarray*}
The left $\mathcal G$-action and right $\mathcal H$-action on this space are defined by
$$(g',g^{-1}\upsilon_0(z))(g,z,h)(h',h'^{-1}h^{-1}\varphi_0(z))=(gg'^{-1},z,hh').$$

We now translate this back to homomorphisms of Lie groupoids, and  
construct the span of homomorphisms corresponding to this bundle, as in \cite{MM}:
$$
\xymatrix{{\mathcal G}&\ar[l] 
	{\mathcal G}\ltimes (R_\upsilon\otimes_{\mathcal K}R_\varphi)\rtimes{\mathcal H} \ar[r] 
		& {\mathcal H}.}
$$
The space of objects in this middle groupoid is 
$$
({\mathcal G}\ltimes (R_\upsilon\otimes_{\mathcal K}R_\varphi)\rtimes{\mathcal H})_0=
R_\upsilon\otimes_{\mathcal K}R_\varphi
$$ 
and the space of arrows is 
$$({\mathcal G}\ltimes (R_\upsilon\otimes_{\mathcal K}R_\varphi)\rtimes{\mathcal H})_1=(G\times X)\times_X (R_\upsilon\otimes_{\mathcal K}R_\varphi)\times_Y(H\times Y) \cong 
G\times (R_\upsilon\otimes_{\mathcal K}R_\varphi)\times H.$$ 
So 
${\mathcal G}\ltimes (R_\upsilon\otimes_{\mathcal K}R_\varphi)\rtimes{\mathcal H}
\cong G\ltimes (R_\upsilon\otimes_{\mathcal K}R_\varphi) \rtimes H
\cong (G\times H)\ltimes (R_\upsilon\otimes_{\mathcal K}R_\varphi)$.
The source map is defined by projection, and the target map is defined by the (left) action of $G\times H$,
$t(g',h',[g,z,h])=[gg'^{-1},z,hh'^{-1}].$

The homomorphisms 
\begin{equation}\label{span2}
\xymatrix@1{
{\mathcal G}&\ar[l]_-\omega (G\times H)\ltimes (R_\upsilon\otimes_{\mathcal K}R_\varphi)\ar[r]^-\psi &{\mathcal H}
}
\end{equation}
are defined by
$$\omega_0[g,z,h]=q_\upsilon(g,z,h)=g^{-1}\upsilon_0(z), \, \, \, \, \, \omega_1(g',h',[g,z,h])=(g',g^{-1}\upsilon_0(z))$$
and $$\psi_0[g,z,h]=q_\varphi(g,z,h)=h^{-1}\varphi_0(z), \, \, \, \, \, \psi_1(g',h',[g,z,h])=(h',h^{-1}\varphi_0(z)).$$

Finally, we construct a 2-cell in the bicategory of fractions from the generalized map in (\ref{span}) to the one in (\ref{span2}).
To this end, define a homomorphism 
$$\theta\colon{\mathcal K}\rightarrow (G\times H)\ltimes (R_\upsilon\otimes_{\mathcal K}R_\varphi)$$
by 
$$\theta_0(z)=[e_G,z,e_H]\mbox{ and }\theta_1(k)=(\pi_1\upsilon_1(k),\pi_1\varphi_1(k),[e_G,s(k),e_H]).$$
We claim that the following diagram of groupoids and homomorphisms commutes:
\begin{equation}\label{2-iso}
\xymatrix{
&{\mathcal K}\ar[dl]_\upsilon\ar[dr]^\varphi\ar[dd]_\theta
\\
G\ltimes X && H\ltimes Y
\\
&(G\times H)\ltimes (R_\upsilon\otimes_{\mathcal K}R_\varphi)\ar[ul]^\omega\ar[ur]_\psi
}
\end{equation}
Indeed, 
\begin{eqnarray*}
\omega_0\circ\theta_0(z)&=&w_0[e_G,z,e_H]\\
&=&\upsilon_0(z),\\
 \omega_1\circ\theta_1(k)&=&w_1( \pi_1\upsilon_1(k),\pi_1\varphi_1(k),[e_G,s(k),e_H] )\\
&=&(\pi_1\upsilon_1(k),\upsilon_0(s(k)))
\\
&=&\upsilon_1(k),
\end{eqnarray*}
and 
\begin{eqnarray*}\psi_0\circ\theta_0(z)&=&\psi_0[e_G,z,e_H]\\&=&\varphi_0(z),\\
\psi_1\circ\theta_1(k)&=&\psi_1(\pi_1\upsilon_1(k),\pi_1\varphi_1(k),[e_G,s(k),e_H])\\
&=&(\pi_1\varphi_1(k),\varphi_0(s(k)))\\
&=&\varphi_1(k).
\end{eqnarray*}
We conclude the proof by remarking that the diagram (\ref{2-iso}) represents an (invertible) 2-cell
in $\mbox{\sf LieGpd}(W^{-1})(G\ltimes X,H\ltimes Y)$.
\end{proof}

The previous proposition implies that for any two orbifold translation groupoids $G\ltimes X$ and $H\ltimes Y$, the inclusion
of categories 
$$
\mbox{\sf EqTrGpd}(W^{-1})(G\ltimes X,H\ltimes Y)\hookrightarrow \mbox{\sf LieGpd}(W^{-1})(G\ltimes X,H\ltimes Y)
$$ 
is essentially surjective on objects, {\it i.e.}, on morphisms $G\ltimes X\rightarrow H\ltimes Y$. 
It remains to be shown in the proof of Theorem \ref{T:main} that the inclusion functor 
$$
\mbox{\sf EqTrGpd}(W^{-1})(G\ltimes X,H\ltimes Y)\hookrightarrow \mbox{\sf LieGpd}(W^{-1})(G\ltimes X,H\ltimes Y)
$$ 
is fully faithful on arrows, {\it i.e.}, on 2-cells between morphisms $G\ltimes X\rightarrow H\ltimes Y$.

\begin{prop}\label{2-cell-correspondence}
Any 2-cell 
$$[{\mathcal M},\theta,\theta',\alpha_1,\alpha_2]\colon((\upsilon,w),K\ltimes Z,(\varphi,f))\Rightarrow((\upsilon',w'),K'\ltimes Z',(\varphi',f'))$$
for orbifold groupoids is equivalent to a 2-cell of the form $[(K\times K')\ltimes Q,\kappa,\kappa',\alpha_1',\alpha'_2]$,
where $\kappa$ and $\kappa'$ are equivariant essential equivalences.
\end{prop}

\begin{proof}
Since $\theta$ is an essential equivalence,  the span 
$K\ltimes Z\stackrel{\theta}{\longleftarrow}{\mathcal M}\stackrel{\theta'}{\longrightarrow}K'\ltimes Z'$
represents a  generalized map from $K\ltimes Z$ to $K'\ltimes Z'$.
We will again use the correspondence with the Hilsum-Skandalis maps  
to find a span of equivariant essential equivalences
which are part of an equivalent 2-cell.
As in the proof of Proposition \ref{morphisms}, we find that 
$R_{\theta}^{-1}\otimes_{\mathcal M}R_{\theta'}\cong (K\times M_0\times K')/\sim_{\mathcal M}$, where the action of
${\mathcal M}$ is defined by 
$$(k,x,k')\cdot m=(\pi_1\theta_1(m^{-1})k,s(m),\pi_1\theta'_1(m^{-1})k'),$$ for $m\in M_1$ with
$t(m)=x\in M_0$. The bundle maps into $Z$ and $Z'$,
$$
\xymatrix{Z & \ar[l]_-{q_{\theta}}\ar[r]^-{q_{\theta'}}R_\theta^{-1}\otimes_{\mathcal K}R_{\theta'} & Z',}
$$
are defined by
$$
q_\theta(k,x,k')=k^{-1}\theta_0(x)\mbox{ and }  q_{\theta'}(k,x,k')=k'^{-1}\theta'_0(x).
$$
The corresponding  span of  equivariant homomorphisms from an intermediate translation groupoid 
into $K\ltimes Z$ and $K'\ltimes Z'$ is given by
$$
\xymatrix@1{K\ltimes Z&\ar[l]_-\kappa (K\times K')\ltimes (R_\theta^{-1}\otimes_{\mathcal M}R_{\theta'})\ar[r]^-{\kappa'} &K'\ltimes Z'}
$$
defined by
$$\kappa_0[k,x,k']=q_\theta(k,x,k')=k^{-1}\theta_0(x),\, \, \,  \kappa_1(\ell,\ell',[k,x,k'])=(\ell,k^{-1}\theta_0(x))$$
and 
$$\kappa'_0[k,x,k']=q_{\theta'}(k,x,k')=k'^{-1}\theta'_0(x), \, \, \, \kappa'_1(\ell,\ell',[k,x,k'])=(\ell',k'^{-1}\theta'_0(x)).$$
So let $Q=R_\theta^{-1}\otimes_{\mathcal M}R_{\theta'}$.

Note that the  natural transformations $\alpha_1$ and $\alpha_2$ are given by smooth functions
$\alpha_1\colon M_0\rightarrow G\times X$ and $\alpha_2\colon M_0\rightarrow H\times Y$.
We will denote the components of these functions by $\alpha_1(x)=(\alpha_1^G(x), \alpha_1^X(x))$ and 
$\alpha_2(x)=(\alpha_2^H(x), \alpha_2^Y(x))$.
We define the new  transformations 
$$\alpha_1'\colon  (R_\theta^{-1}\otimes_{\mathcal M}R_{\theta'})\rightarrow G\times X,\quad 
\alpha'_2\colon (R_\theta^{-1}\otimes_{\mathcal M}R_{\theta'})\rightarrow H\times Y$$
by
$$\alpha_1'[k,x,k']=(\upsilon'(k')^{-1}\alpha_1^G(x)\upsilon(k),\upsilon(k)^{-1} w(\theta_0(x))),$$
and
$$\alpha_2'[k,x,k']=(\varphi'(k')^{-1}\alpha_2^H(x)\varphi(k),\varphi(k)^{-1} f(\theta_0(x))).$$
The fact that $\alpha_1'$ and $\alpha_2'$ are well-defined on equivalence classes follows from the
fact that $\alpha_1$ and $\alpha_2$ satisfy the naturality condition, as the following
calculation shows,
\begin{eqnarray*}
\lefteqn{\alpha_1' (\pi_1\theta (m^{-1})k,s(m),\pi_1\theta' (m^{-1})k')=}\\
&=&
(\upsilon' (\pi_1\theta' (m^{-1})k')^{-1}\alpha_1^G (s(m))\upsilon (\pi_1\theta ( m^{-1} ) k),
\upsilon (\pi_1\theta (m^{-1})k)^{-1} w(\theta_0 (s(m))))
\\
&=&
(\upsilon'(k')^{-1} \upsilon' ( \pi_1\theta'_1(m) )\alpha_1^G ( s(m))\upsilon (\pi_1\theta (m^{-1}))
\upsilon(k),\upsilon (k)^{-1}\upsilon (\pi_1\theta_1 (m)) w(\theta_0 (s(m))))
\\
&=&
(\upsilon' (k')^{-1}\alpha_1^G (t(m))\upsilon (k),\upsilon (k)^{-1} w(\theta_0 (t(m))))
\\
&=&
(\upsilon' (k')^{-1}\alpha_1^G (x)\upsilon (k),\upsilon (k)^{-1} w(\theta_0 (x)))
\\
&=&
\alpha_1'(k,x,k').
\end{eqnarray*}
The fact that $\alpha_1'$ and $\alpha_2'$ satisfy the naturality condition 
can be checked by a straightforward calculation.
Also, 
$$s\circ\alpha_1'  [k,x,k']=\upsilon(k)^{-1}\cdot w(\theta_0(x))=w(\kappa_0(x)),$$ 
and 
\begin{eqnarray*}
t\circ\alpha_1'[k,x,k']&=&\upsilon'(k')^{-1}\alpha_1^G(x)\upsilon(k) (\upsilon(k)^{-1} w(\theta_0(x)))\\
&=&
\upsilon'(k')^{-1}\alpha_1^G(x) w(\theta_0(x))\\
&=&\upsilon'(k')^{-1} w'(\theta'_0(x))\\
&=& w'(k'^{-1} \theta'_0(x))\\
&=&w'\kappa'_0(x),
\end{eqnarray*}
so $\alpha_1'$ represents a natural transformation from 
$(\upsilon,w)\circ\kappa$ to $(\upsilon',w')\circ\kappa'$.
The calculation for $\alpha_2'$ goes similarly.
\end{proof}

\begin{rmk}
We have only shown that the inclusion functor $\mbox{\sf EqTrGpd}(W^{-1})\hookrightarrow
\mbox{\sf LieGpd}_{\mbox{\scriptsize\sf Tr}}(W^{-1})$ is a (weak) equivalence of bicategories,
and this is sufficient for our purposes. However, the method of the proof can also be used to 
construct a homomorphism of bicategories 
$\Phi\colon\mbox{\sf LieGpd}_{\mbox{\scriptsize\sf Tr}}(W^{-1})\rightarrow\mbox{\sf EqTrGpd}(W^{-1})$
in the opposite direction. On objects, $\Phi$ is the identity, and it sends a generalized morphism
$$
\xymatrix{
G\ltimes X & {\mathcal K}\ar[l]_-{\upsilon}\ar[r]^-\varphi & H\ltimes Y
}
$$ 
to 
$$
\xymatrix{
G\ltimes X & (G\times H)\ltimes(R_{\upsilon}^{-1}\otimes_{\mathcal K}R_\varphi) 
		\ar[l]_-{(\pi_1, q_\upsilon)}\ar[r]^-{(\pi_2, q_\varphi)} & H\ltimes Y\rlap{,}
}
$$ 
as constructed above. For a 2-cell 
\begin{equation}\label{new2-cell}
\xymatrix@R=2em@C=3em{
	 &{\mathcal K}\ar[dl]_{\omega}\ar[dr]^\varphi
 \\
 G\ltimes X \ar@{}[r]|{\alpha_1\Downarrow} 
	& {\mathcal L}\ar[u]_{\nu}\ar[d]^{\nu'}\ar@{}[r]|{\alpha_2\Downarrow} & H\ltimes Y
 \\
	 &{\mathcal K}'\ar[ul]^{\omega'}\ar[ur]_{\varphi'} &,
 }
\end{equation}
consider the induced 2-cell
\begin{equation}\label{fat2-cell}
\xymatrix@R=2em@C=3em{
	& (G\times H)\ltimes(R_{\omega}^{-1}\otimes_{\mathcal K}R_\varphi) \ar[ddl]_{(\pi_1, q_\omega)}
		\ar[ddr]^{(\pi_2, q_\varphi)}&
\\
	 &{\mathcal K}\ar[u]_{\theta} \ar[dl]_{\omega}\ar[dr]^\varphi \ar@{}[ul]|(.2){=}\ar@{}[ur]|(.2){=}
 \\
 G\ltimes X \ar@{}[r]|{\alpha_1\Downarrow} 
	& {\mathcal L}\ar[u]_{\nu}\ar[d]^{\nu'}\ar@{}[r]|{\alpha_2\Downarrow} & H\ltimes Y
 \\
	 &{\mathcal K}'\ar[ul]^{\omega'}\ar[ur]_{\varphi'} \ar[d]^{\theta'} \ar@{}[dl]|(.2){=}\ar@{}[dr]|(.2){=}&
\\
	& (G\times H)\ltimes(R_{\omega'}^{-1}\otimes_{\mathcal K'}R_{\varphi'}) \ar[uul]^{(\pi_1, q_{\omega'})}
		\ar[uur]_{(\pi_2, q_{\varphi'})} &,
 }
\end{equation}
where $\theta$ and $\theta'$ are the morphisms as described in (\ref{2-iso}).
Then $\Phi$ sends (\ref{new2-cell}) to the 2-cell
$$
\xymatrix@R=2em@C=3em{
	 &(G\times H)\ltimes(R_{\upsilon}^{-1}\otimes_{\mathcal K}R_\varphi) 
		\ar[dl]_{(\pi_1,q_{\omega})}\ar[dr]^{(\pi_2, q_{\varphi'})}
 \\
 G\ltimes X \ar@{}[r]|{\alpha'_1\Downarrow} 
	& K\ltimes Z
		\ar[u]_{\kappa}\ar[d]^{\kappa'}\ar@{}[r]|{\alpha'_2\Downarrow} & H\ltimes Y
 \\
	 &(G\times H)\ltimes(R_{\omega'}^{-1}\otimes_{\mathcal K'}R_{\varphi'}) 
		\ar[ul]^{(\pi_1,q_{\omega'})}\ar[ur]_{(\pi_2, q_{\varphi'})} &,
 }
$$
obtained by applying the methods of the proof of Proposition \ref{2-cell-correspondence} to
(\ref{fat2-cell}).
\end{rmk}

\section{Proofs III:  Proofs of Additional Results} \label{S:morepfs}

In this section we include proofs of the additional results mentioned throughout the paper.  
We begin with the lemma that the fibre product of 
two translation groupoids is another translation groupoid.  

\begin{proof}[Proof of Lemma \ref{L:pb}]
The object space of the fibre product groupoid (if it exists) is $P_0= Y\times_{K_0}K_1\times_{K_0} X$,
so its elements can be represented by triples 
$$
(y,\varphi_0(y)\stackrel{k}{\rightarrow}\psi_0(x),x),
$$ 
where 
$y\in Y$, $k\in K_1$, and $x\in X$.
An element of the space of arrows $P_1$ is given by a triple 
$$
(y\stackrel{(h,y)}{\longrightarrow}hy,\varphi_0(y)\stackrel{k}{\longrightarrow}\psi_0(x),x\stackrel{(g,x)}{\longrightarrow}gx)
$$
with $y\in Y$, $h\in H$, $k\in K$, $x\in X$ and $g\in G$.
Such triples are in one-to-one correspondence with 5-tuples of the form
$(h,y,\varphi_0(y)\stackrel{k}{\rightarrow}\psi_0(x),x,g)$. Moreover, in this notation, 
$$
s(h,y,\varphi_0(y)\stackrel{k}{\rightarrow}\psi_0(x),x,g)=(y,\varphi_0(y)\stackrel{k}{\rightarrow}\psi_0(x),x),
$$
and 
$$
t(h,y,q_0(y)\stackrel{k}{\rightarrow}\psi_0(x),x,g)=(hy,\varphi_0(hy)\stackrel{\psi_1(g,x)k[\varphi_1(h,y)]^{-1}}{\longrightarrow} \psi_0(gx),gx),
$$
so ${\mathcal P}$ is the translation groupoid for the action of $G\times H$ on $P_0=Y\times_{K_0}K_1\times_{K_0} X$,
defined by
$$
(g,h)\cdot(y,\varphi_0(y)\stackrel{k}{\rightarrow}\psi_0(x),x)=
	(hy,\varphi_0(hy)\stackrel{\psi_1(g,x)k[\varphi_1(h,y)]^{-1}}{\longrightarrow} \psi_0(gx),gx).
$$
Also, $\xi_0(y,k,x)=y$, $\xi_1(h,y,k,x,g)=(h,y)$, $\zeta_0(y,k,x)=x$, and $\zeta_1(h,y,k,x,g)=(g,x)$, so these maps have the desired format.
\end{proof}

Next we prove that all equivariant essential equivalences between translation groupoids have the forms  specified in Proposition \ref{L:form}.

\begin{proof}[Proof of Proposition \ref{L:form}]    Let 
$$
\xymatrix{
G\times X \ar@<-.4ex>[d] \ar@<.4ex>[d] \ar[r]^{\varphi\times f} & H\times Y \ar@<-.4ex>[d] \ar@<.4ex>[d]
\\
X\ar[r]_f & Y
}
$$
be an equivariant essential equivalence between translation groupoids. 
We will denote this by $\varphi\ltimes f\colon G\ltimes X\rightarrow H\ltimes Y$.
This map  can be factored in the following way:
$$
\xymatrix@C=3.5em{
G\times X \ar@<-.4ex>[d] \ar@<.4ex>[d] \ar[r]^-{\overline\varphi\times f} & 
	G/\mbox{Ker}(\varphi)\times f(X) \ar@<-.4ex>[d] \ar@<.4ex>[d] \ar[r]^-{\mbox{\scriptsize inclusion}} & 
	H\times Y \ar@<-.4ex>[d] \ar@<.4ex>[d]
\\
X\ar[r]_f & f(X) \ar[r]_-{\mbox{\scriptsize inclusion}}& Y.
}
$$

Since the map $\overline\varphi\ltimes f$ is surjective on objects and $\varphi\ltimes f$ is essentially surjective,
so is $\overline\varphi\ltimes f$. Similarly, the  right inclusion map is essentially surjective because   $\varphi\ltimes f$ is.
  
  We will show that 
with the notation above,  the first map $\overline\varphi\ltimes f$ is of the form

$$ (\ref{quotientform}) \phantom{spacespace} G\ltimes X \rightarrow G/K\ltimes X/K,$$ 
where $K$ is a normal subgroup of $G$ which acts freely on $X$, and $X/K$ is the
quotient of $X$ by this action. The second map is of the form 
$$(\ref{inclform})  \phantom{spacespace}
K\ltimes Z\rightarrow H\ltimes (H\times_K Z),
$$
where $K$ is a 
(not necessarily normal) subgroup of $H$.

Consider the diagram
$$
\xymatrix{
G\times X \ar[d]_{(s,t)} \ar[r]^-{\overline\varphi\times f} & G/\mbox{Ker}(\varphi)\times f(X)\ar[d]_{(s,t)} \ar[r]&
	H\times Y \ar[d]_{(s,t)}
\\
X\times X\ar[r]_-{f\times f} & f(X)\times f(X)\ar[r]_-{\mbox{\scriptsize incl}\times\mbox{\scriptsize incl}} & Y\times Y.
}
$$

We show  that the right hand square is a pullback. 
Let $p\colon P\rightarrow H\times Y$ and $q\colon P\rightarrow f(X)\times f(X)$
be such that $(s,t)\circ p=(\mbox{incl}\times\mbox{incl})\circ q$. 
Then there is a map $r\colon P\rightarrow G/\mbox{Ker}(\varphi)\times f(X)$
defined as follows: let $\pi\in P$, and let $p(\pi)=(h_\pi,y_\pi)$ and $q(\pi)=(y'_\pi,y''_\pi)$. 
Then $y'_\pi=y_\pi$ and $y''_\pi=h_\pi y_\pi$.
Choose $x$ and $x'$ in $X$ such that $f(x)=y'_\pi$ and $f(x')=y''_\pi$ . 
Since $\varphi\ltimes f$ is a essential equivalence, 
there is a unique $g\in G$ such that $gx=x'$ and $\varphi(g)=h$. 
We define $r(\pi)=(\overline{g}, y'_\pi)$.
To show that this does not depend on the choice of the pre-images $x$ and $x'$, let $z$ and $z'$ 
be such that $f(z)=y'_\pi$ and $f(z')=y''_\pi$,
and let $g'\in G$ be the unique element such that $g'z=z'$ and $\varphi(g')=h_\pi$.
Since $f(z)=f(x)$ and $f(z')=f(x')$, and $\varphi\ltimes f$ is a essential equivalence, there are unique elements 
$a,a'\in G$ such that $ax=z$, $a'x'=z'$ and $\varphi(a)=e_H=\varphi(a')$. 
Moreover, $g'a=a'g$, since $g'ax=z'$ and $a'gx=z'$, 
and $\varphi(g'a)=h_\pi=\varphi(a'g)$. So $\overline{g}=\overline{g'}\in G/\mbox{Ker}(\varphi)$.
It is clear that the map $r\colon P\rightarrow G/\mbox{Ker}(\varphi)\times f(X)$ is the unique map which makes the
following diagram commute:
$$
\xymatrix{
P\ar@/_/[ddr]_{q}\ar@/^/[drr]^{p}\ar[dr]^r 
\\
& G/\mbox{Ker}(\varphi)\times f(X) \ar[d]\ar[r] & H\times Y\ar[d]
\\
&f(X)\times f(X) \ar[r] & Y\times Y,
}
$$ 
so the square is a pullback. Since $\varphi\ltimes f$ is an essential equivalence, 
the whole rectangle is also a pullback, so the left hand
square is a pullback. 

We conclude that we have factored $\varphi\ltimes f$ into two new essential equivalences.  It is easy to check that $\overline{\varphi} \ltimes f$ has the form of a projection $G \ltimes X \to G/\mbox{Ker}(\varphi) \ltimes X/ \mbox{Ker}(\varphi) $.  So it remains to show that the space  
$Y$ is homeomorphic to the group extension of the $G/\mbox{Ker}(\varphi)$-space $f(X)$ over the 
inclusion $G/\mbox{Ker}(\varphi) \to H$, that is, that $Y\cong H\times_{G/\mbox{Ker}(\varphi)} f(X)$.
 
Note that elements of $H\times_{G/\mbox{Ker}(\varphi)} f(X)$ are represented by pairs $(h,f(x))$  with $h\in H$ and $x\in X$,
and $(h\varphi(\overline{g}),f(x))\sim (h,\varphi(\overline{g})f(x))$.
There is a morphism $H\times_{G/\mbox{Ker}(\varphi)} f(X)\rightarrow Y$, defined by $(h,f(x))\mapsto hf(x).$
This map is a surjective submersion since $\varphi\ltimes f$ is essentially surjective. It is also injective:
if $hf(x)=h'f(x')$, then $h'^{-1}hf(x)=f(x')$, so there is an element $g\in G$ such that $gx=x'$ and $\varphi(g)=h'^{-1}h$,
so $h=h'\varphi(g)$. So $H\times_{G/\mbox{Ker}(\varphi)} f(X)\cong Y$, as desired.

We conclude that all essential equivalences can be obtained as composites of essential equivalences of the forms (\ref{quotientform})
and (\ref{inclform}).  \end{proof}

Finally, we include the proof of the 2-for-3 Lemma mentioned in Section \ref{S:back}.  
\begin{lem}\label{2outof3}
The class of essential equivalences between Lie groupoids satisfies the 2-for-3 property, {\it i.e.}, if we have 
 homomorphisms
$
{\mathcal G}\overset{\varphi}\to  {\mathcal K}\overset{\psi}{\to} {\mathcal H}$
 such that 
two out of $\{\varphi,\psi,\varphi\circ\psi\}$ are essential equivalences, then so is the third.
\end{lem}

\begin{proof}
Consider the following diagram
$$
\xymatrix{
G_1\ar[d]_{(s,t)}\ar[r]^{\varphi_1} \ar@{}[dr]|{(A)}& K_1 \ar[d]_{(s,t)}\ar[r]^{\psi_1} \ar@{}[dr]|{(B)} & H_1 \ar[d]^{(s,t)}
\\
G_0\times G_0 \ar[r]_{\varphi_0\times\varphi_0} & K_0\times K_0 \ar[r]_{\psi_0\times\psi_0} & H_0\times H_0.
}
$$
It is a standard property of fibre products that if any two out of $(A)$, $(B)$, and the whole square are fibre products, so is the third.
So if any two out of
$\{\varphi,\psi,\varphi\circ\psi\}$ are fully faithful, then so is the third.

It is straightforward to show that if $\varphi$ and $\psi$ are essentially surjective, so is the composite $\psi\circ\varphi$.
It is also straightforward to show that if $\psi\circ\varphi$ is essentially surjective, then $\psi$ is essentially surjective.

Lastly, suppose  that $\psi$ and $\varphi\circ\psi$ are essential equivalences. We claim that this implies that $\varphi$
is essentially surjective (and therefore an essential equivalence). Since $\psi$ is fully faithful, we have that 
$K_1\cong K_0\times_{H_0,s}H_1\times_{t,H_0}K_0$, and therefore 
$G_0\times_{\varphi_0,K_0,s}K_1\cong G_0\times_{K_0}K_0\times_{H_0}H_1\times_{H_0}K_0\cong G_0\times_{H_0}H_1\times_{H_0}K_0$.
So consider the following commutative diagram.
$$
\xymatrix{
G_0\times_{K_0}K_1 \ar[d]_{\wr} \ar[r]^{\pi_2} & K_1 \ar[dr]^t
\\
G_0\times_{H_0}H_1\times_{H_0}K_0 \ar[rr]^{\pi_3} \ar@{}[dr]|{\mbox{\scriptsize pb}} \ar[d] && K_0\ar[d]^{\psi_0}
\\
G_0\times_{H_0}H_1 \ar[d]_{\pi_1} \ar[r]_{\pi_2} \ar@{}[dr]|{\mbox{\scriptsize pb}}  & H_1\ar[r]_t \ar[d]^s & H_0
\\
G_0\ar[r]_{\psi_0\circ\varphi_0} & H_0
}
$$
The composite $t\circ\pi_2\colon G_0\times_{H_0}H_1\rightarrow H_0$ is a surjective submersion, because $\varphi\circ\psi$ is essentially surjective.
So, $\pi_3\colon G_0\times_{H_0}H_1\times_{H_0}K_0\rightarrow K_0$ is a surjective submersion, since it is a pullback of one,
and this makes $t\circ\pi_2\colon.G_0\times_{K_0}K_1\rightarrow K_0$ a surjective submersion.
We conclude that  in this case $\varphi$ is also essentially surjective.
\end{proof}

\end{document}